\documentclass{elsart}

\usepackage{amssymb}
\usepackage{epsfig}
\usepackage{amsmath}
\usepackage{amsfonts}
\newtheorem{theorem}{Theorem}

\newtheorem{lemme}{Lemma}

\begin{document}

\begin{frontmatter}

\title{Modeling 1-D elastic P-waves in a fractured \\rock with hyperbolic jump conditions}

\author{Bruno Lombard},
\ead{lombard@lma.cnrs-mrs.fr}
\ead[url]{http://w3lma.cnrs-mrs.fr/$\sim$MI/}
\author{Jo\"el Piraux}
\ead{piraux@lma.cnrs-mrs.fr}
\address{Laboratoire de M\'ecanique et d'Acoustique,\\ 31 chemin Joseph Aiguier, 13402 Marseille, France}

\begin{abstract}
The propagation of elastic waves in a fractured rock is investigated, both theoretically and numerically. Outside the fractures, the propagation of compressional waves is described in the simple framework of one-dimensional linear elastodynamics. The focus here is on the interactions between the waves and fractures: for this purpose, the mechanical behavior of the fractures is modeled using nonlinear jump conditions deduced from the Bandis-Barton model  classicaly used in geomechanics. Well-posedness of the initial-boundary value problem thus obtained is proved. Numerical modeling is performed by coupling a time-domain finite-difference scheme with an interface method accounting for the jump conditions. The numerical experiments show the effects of contact nonlinearities. The harmonics generated may provide a non-destructive means of evaluating the mechanical properties of fractures.
\end{abstract}
\begin{keyword}
elastic waves, contact nonlinearity, Bandis-Barton model, jump conditions, finite-difference schemes, interface method.
\PACS 02.60.Cb, 02.70.Bf, 43.25.+y, 46.50.+a
\end{keyword}
\end{frontmatter}

\section{Introduction}\label{SecIntro}

Fractures are the breaks in rocks caused by the huge stresses resulting from plate tectonics. It is of fundamental importance for geophysicists to be able to determine the position and the properties  of fractures (such as their thickness) to be able to make predictions about the mechanical properties of a fractured platform or the diffusion of a pollutant, for instance. Elastic waves are commonly used for this purpose. When the wavelengths are much larger than the distance between fractures, the fractures are generally not studied individually, and homogenization theories are applied. Otherwise, as in the case of the present study, it is possible to study single fractures. If, in addition, the wavelengths are much larger than the thickness of the fractures, the latter can be modeled in terms of interfaces with appropriate jump conditions.

Many experimental, theoretical and numerical studies have dealt with wave propagation across thin and single fractures in terms of linear jump conditions \cite{PYRAK90,ROKHLIN1,ALIMENTAIRE1}. The linear framework provides an appealing approach but it may not be very realistic, since non-physical penetration of both sides of the fractures may occur. Some authors have proposed more accurate fracture models, such as the Bandis-Barton model \cite{BANDIS83}. This model is commonly used in rock mechanics and engineering to deal with quasistatic loading conditions, such as those occuring in flood barriers, for instance. However, very few studies have dealt so far with wave propagation in a model of this kind: we have only found one study using this appoach in \cite{ZHAO01}. 

The aim of the present paper is to further study the propagation of mechanical waves across a fracture described by the Bandis-Barton model, both theoretically and numerically, in a highly idealized configuration. To perform the numerical modeling, we combine an {\it interface method} with a classical numerical scheme for wave propagation, namely the ADER scheme \cite{ADER04}. The interface method involves changing this scheme locally, using the jump conditions at the interface; moreover, it gives a subcell resolution when the interface does not coincide with the meshing. Many interface methods have been proposed since the 90's; see \cite{LI_IIM_03} for a review. Here, we adapt the explicit simplified interface method (ESIM) previously developed for dealing with linear contacts \cite{ALIMENTAIRE1}; see \cite{BIBLE1} for an overview of the principles and advantages of this method. To our knowledge, this is the first time nonlinear jump conditions have been studied using an interface method

The paper is organized as follows. In section 2, the problem under study is introduced, especially the hyperbolic jump conditions (\ref{JCBB}). An analysis of the solution is performed in section \ref{SecMath}: conservation of energy, and the existence, uniqueness and regularity of the solution. Section  \ref{SecNum} deals with the numerical methods. The numerical experiments performed in section \ref{SecExp} with realistic parameters show the influence of the incident wave amplitudes. 

\section{Problem statement}\label{SecPS}

\subsection{Configuration}\label{SecConfig}

\begin{figure}[htbp]
\begin{center}
\includegraphics[scale=0.9]{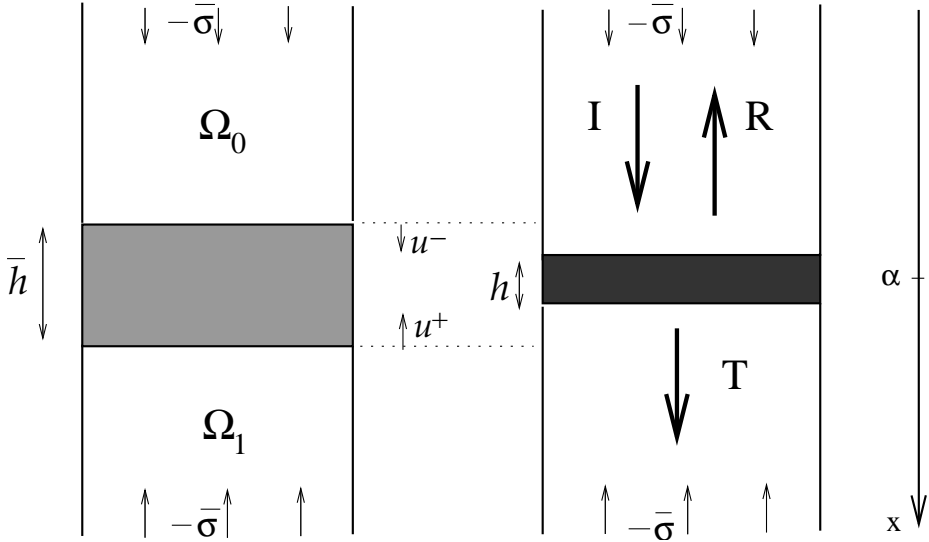}
\caption{Static (left) and dynamic (right) behavior of the fractured rock. I: incident wave; R: reflected wave; T: transmitted wave.}
\label{FigStatDyn}
\end{center}
\end{figure}

Consider a rock with a single plane fracture. Outside the fracture, the media involved $\Omega_i$ ($i=0,1$) are linearly elastic and isotropic; they are subject to a constant static stress $-\overline{\sigma}$ $(\overline{\sigma}>0)$ running perpendicular to the fracture. At rest, the fractured zone is an {\it interphase} with thickness $\overline{h}>0$ (figure \ref{FigStatDyn}, left). The nonlinear mechanical behavior of the interphase is investigated in section \ref{SecJC0}.

Let us now consider a plane compressional wave propagating through $\Omega_0$ normal to the interphase; the interactions between this incident wave and the interphase give rise to reflected (in $\Omega_0$) and transmitted (in $\Omega_1$) plane compressional waves. These perturbations in $\Omega_0$ and $\Omega_1$ are described by the simple one-dimensional elastodynamic equations \cite{ACHENBACH}
\begin{equation}
\rho \,\frac{\textstyle \partial \,v}{\textstyle \partial \,t}=\frac{\textstyle \partial\,\sigma}{\textstyle \,\partial\, x},\qquad
\frac{\textstyle \partial \,\sigma}{\textstyle \partial \,t}=
\rho \,c^2 \,\frac{\textstyle \partial \,v}{\textstyle \partial \,x},
\label{LCscal}
\end{equation}
where $v=\frac{\partial\,u}{\partial\,t}$ is the elastic velocity, $u$ is the elastic displacement, and $\sigma$ is the elastic stress perturbation around $-\overline{\sigma}$. The physical parameters involved are the density $\rho$ and the elastic speed of the compressional waves $c$; these piecewise constant parameters may be discontinuous around the fracture: $(\rho_0,c_0)$ if $x \in \Omega_0$, $(\rho_1,c_1)$ if $x \in \Omega_1.$ The dynamic stresses induced by the elastic waves affect the thickness $h(t)$ of the interphase (figure \ref{FigStatDyn}, right). Due to the finite compressibility of the interphase, the constraint
\begin{equation}
h=\overline{h}+[u]\geq \overline{h}-d>0
\label{Penetration}
\end{equation}
must be satisfied, where $[u]=u^+-u^-$ is the difference between the elastic displacements on the two sides of the interphase, and $d>0$ is the {\it maximum allowable closure} \cite{BANDIS83}. We also assume that the wavelengths of the elastic perturbations are much larger than $h$. One can therefore neglect the propagation time through the interphase, and replace it by a zero-thickness {\it interface} at $x=\alpha$, where $\alpha$ belongs to the interphase; therefore, $[u]=[u(\alpha,\,t)]=u(\alpha^+,\,t)-u(\alpha^-,\,t)$. 

\subsection{Bandis-Barton model}\label{SecJC0}

\begin{figure}[htbp]
\begin{center}
\includegraphics[scale=0.9]{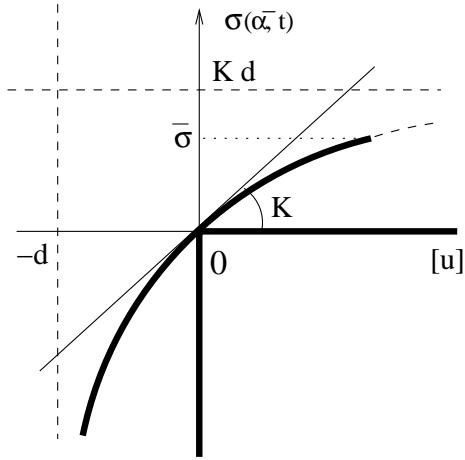}
\end{center}
\caption{Sketch of the stress-displacement relation deduced from the Bandis-Barton model. The bold straight segments denote the limit-case of unilateral contact.}
\label{FigBB}
\end{figure}

Single thin fractures have been classically modeled in terms of linear jump conditions \cite{PYRAK90}. Given a {\it stiffness} $K>0$ and neglecting the inertial effects, the most usual linear jump conditions are
\begin{equation}
\begin{array}{l}
\left[\sigma(\alpha,\,t)\right]= 0,\\
[8pt]
\displaystyle
\left[u(\alpha,\,t)\right]=\frac{\textstyle 1}{\textstyle K}\,\sigma(\alpha^-,\,t).
\end{array}
\label{JClin}
\end{equation}
The simple jump conditions (\ref{JClin}) can be rigorously obtained by performing an asymptotic analysis of the wave propagation process within a plane interphase which is much thinner than the wavelength ($h \ll \lambda$); then $K=\rho\,c^2/h$, where $\rho$ and $c$ are the physical parameters of the interphase. For $K \rightarrow +\infty$, we obtain perfectly-bonded conditions; for $K \rightarrow 0^+$, we obtain $\sigma(\alpha^\pm,\,t)\rightarrow 0$, and hence the two media $\Omega_0$ and $\Omega_1$ tend to be disconnected. The main drawback of the conditions (\ref{JClin}) is that they do not satisfy (\ref{Penetration}) under large compression loadings: $\sigma(\alpha^-,\,t) < - K\,d \, \Rightarrow \, h <\overline{h}-d$, which contradicts (\ref{Penetration}). Hence, (\ref{JClin}) is realistic only in the case of very small perturbations. When larger ones are involved, nonlinear jump conditions are required.

To satisfy (\ref{Penetration}), we use the Bandis-Barton model \cite{BANDIS83}. This model is based on quasi-static compressional experiments showing that the closure of a fracture depends hyperbolically on the stress applied. In the case of dynamic problems \cite{ZHAO01}, the hyperbolic jump conditions can be written 
\begin{equation}
\begin{array}{l}
\left[\sigma(\alpha,\,t)\right]=0,\\
[4pt]
\displaystyle
\left[u(\alpha,\,t)\right]= \frac{\textstyle 1}{\textstyle K}\,\frac{\textstyle \sigma(\alpha^-,\,t)}{\textstyle \displaystyle 1-\frac{\textstyle \sigma(\alpha^-,\,t)}{\textstyle K\,d}},
\end{array}
\label{JCBB}
\end{equation}
with $\sigma(\alpha^-,\,t)<K\,d$ (the other branch of the hyperbola is not realistic). In the Bandis-Barton model \cite{BANDIS83}, the parameters $K$ and $d$ are linked to $\overline{\sigma}$, satisfying $\overline{\sigma}<K\,d$. The second relation in (\ref{JCBB}) is sketched in figure \ref{FigBB}. Under compression loadings: $\sigma(\alpha^\pm,\,t)<0$, the second equation of (\ref{JCBB}) implies: $[u(\alpha,\,t)]>-d$, hence (\ref{Penetration}) is satisfied. Under traction loadings: $\sigma(\alpha^\pm,\,t)>0$, (\ref{Penetration}) is trivially satisfied. In the latter case, nothing prevents $\sigma(\alpha^-,\,t)\geq \overline{\sigma}$, leading to disconnection or adhesion processes; more realistic models that account for such processes are not investigated here. The straight line with a slope $K$ tangential to the hyperbola at the origin describes the linear jump conditions (\ref{JClin}); as deduced from (\ref{JCBB}), the linear conditions are valid only if $|\sigma(\alpha^\pm,\,t)| \ll K\,d$. In the limit-case $d\rightarrow 0^+$ with $K$ bounded,   the hyperbola tends towards the nondifferentiable graph of the unilateral contact, denoted by bold straight segments in figure \ref{FigBB}. The latter graph corresponds to the Signorini's conditions \cite{THESE_SCAROLE}: $\sigma(\alpha^+,\,t)=\sigma(\alpha^-,\,t)\leq 0$, $\left[u(\alpha,\,t)\right]\geq 0$, and $\sigma(\alpha^\pm,\,t)\,\left[u(\alpha,\,t)\right]=0$. This limit-case, which leads to a difficult mathematical analysis and requires suitable numerical tools, is not investigated here. 

\subsection{Initial-boundary value problem}\label{SecIBVP}

From now, we follow a velocity-stress formulation of elastodynamics. It remains therefore to know the jump condition satisfied by $v$ at $\alpha$. For that purpose, we differentiate the second equation of (\ref{JCBB}) with respect to $t$ 
\begin{equation}
v(\alpha^+,\,t) = v(\alpha^-,\,t)+ \frac{\textstyle 1}{\textstyle K}\,\frac{\textstyle 1}{\textstyle \left(\displaystyle 1-\frac{\textstyle \sigma(\alpha^-,\,t)}{\textstyle K\,d}\right)^2}\,\frac{\textstyle \partial\,\sigma}{\textstyle \partial \,t}(\alpha^-,\,t),
\label{JCUt}
\end{equation}
and we replace the time derivative in (\ref{JCUt}) by a spatial derivative via (\ref{LCscal}). Together with the first equation of (\ref{JCBB}), it gives 
\begin{equation}
\begin{array}{l}
\displaystyle
v(\alpha^+,\,t) = \displaystyle v(\alpha^-,\,t)+\frac{\textstyle \rho_0\,c_0^2}{\textstyle K}\,\frac{\textstyle 1}{\textstyle \left(\displaystyle 1-\frac{\textstyle \sigma(\alpha^-,\,t)}{\textstyle K\,d}\right)^2}\,\frac{\textstyle \partial\,v}{\textstyle \partial\,x}(\alpha^-,\,t),\\
\sigma(\alpha^+,\,t)=\sigma(\alpha^-,\,t).
\end{array}
\label{JCBBbis}
\end{equation}
To sum up the conservation laws (\ref{LCscal}) together with the jump conditions (\ref{JCBBbis}), we set up
\begin{equation}
\boldsymbol{U}(x,\,t)=\left(
\begin{array}{c}
v\\
[4pt]
\sigma
\end{array}
\right),
\qquad
\boldsymbol{A}(x)=\left(
\begin{array}{cc}
0           &  \displaystyle -\frac{\textstyle 1}{\textstyle \rho}\\
[4pt]
-\rho \,c^2 & 0
\end{array}
\right),
\label{MatA}
\end{equation}
and then we state the following initial-boundary value problem
\begin{equation}
\left\{
\begin{array}{l}
\displaystyle
\frac{\textstyle \partial}{\textstyle \partial \,t}\,\boldsymbol{U} + \boldsymbol{A}\, \frac{\textstyle \partial}{\textstyle \partial \,x}\,\boldsymbol{U}=\boldsymbol{0}\quad \mbox{ for } x\in \mathbb{R},\quad x\neq \alpha, \quad t\geq t_0,\\
[8pt]
\displaystyle
\boldsymbol{U}(\alpha^+,\,t)=\boldsymbol{D}_0 \left(\boldsymbol{U}(\alpha^-,\,t),\,\frac{\textstyle \partial}{\textstyle \partial\,x}\boldsymbol{U}(\alpha^-,\,t)\right), \\
[8pt]
\displaystyle
\boldsymbol{U}(x,\,t_0)=\boldsymbol{U}_0(x)
\quad \mbox{ for } x \in \mathbb{R},
\end{array}
\right.
\label{IBVP}
\end{equation}
where $\boldsymbol{D}_0:\mathbb{R}^4\rightarrow\mathbb{R}^2$ is a nonlinear application deduced from (\ref{JCBBbis}) and (\ref{MatA}). We assume that the initial data $\boldsymbol{U}_0(x):\mathbb{R}\rightarrow\mathbb{R}^2$ is a $C^p_c$ function with a compact support included in $\Omega_0$, with $p\geq 1$. 

For use in section \ref{SecESIM}, we need also the jump conditions satisfied by some of the spatial derivatives of $\boldsymbol{U}$ (under suitable assumptions of regularity, as defined in theorem \ref{ThExistence}). Differentiating (\ref{JCBBbis}) $m-1$ times with respect to time, and then replacing the time derivatives by spatial derivatives via (\ref{LCscal}), yield nonlinear $m$-th order jump conditions that can be written
\begin{equation}
\frac{\textstyle \partial^m}{\textstyle \partial\,x^m}\boldsymbol{U}(\alpha^+,\,t)=\boldsymbol{D}_m \left(\boldsymbol{U}(\alpha^-, \,t),...,\,\frac{\textstyle \partial^m}{\textstyle \partial\,x^m}\boldsymbol{U}(\alpha^-, \,t), \,\frac{\textstyle \partial^{m+1}}{\textstyle \partial\,x^{m+1}}\boldsymbol{U}(\alpha^-, \,t)\right),
\label{JCBBm}
\end{equation}
where $\boldsymbol{D}_m:\mathbb{R}^{2(m+2)}\rightarrow\mathbb{R}^2$ is a nonlinear application. The computation of $\boldsymbol{D}_m$ is tedious task that can be automated using computer algebra tools. The simulations shown in section \ref{SecExp} were obtained in this way. 

\section{Analysis of the solution}\label{SecMath}

The aim of this subsection is to prove that the initial-boundary value problem (\ref{IBVP}) is a well-posed problem. Our proof is based on rather elementary tools, and yields an analytical solution. The numerical evaluation of the latter with a fourth-order Runge-Kutta scheme constitutes the semi-analytical solution used in section \ref{SecExp}. Before giving the theorem, we need some intermediate results. 

\subsection{Conservation of energy}

In the first lemma, we define an energy, and we show that it is conserved.
\begin{lemme}
Let $\boldsymbol{U}(x,\,t)$ be a solution of (\ref{IBVP}). Then, 
\begin{equation}
\begin{array}{l}
\displaystyle
E(\boldsymbol{U},\,t)=\frac{\textstyle 1}{\textstyle 2}\int_{-\infty}^{+\infty}\left(\rho\,v^2+\frac{\textstyle 1}{\textstyle \rho\,c^2}\,\sigma^2\right)dx\\
[4pt]
\displaystyle
\qquad \qquad +K\,d^2\left(\ln\left(1-\frac{\textstyle \sigma(\alpha^-,\,t)}{\textstyle K\,d}\right)+\frac{\textstyle 1}{\textstyle 1 - \displaystyle \frac{\textstyle \sigma(\alpha^-,\,t)}{\textstyle K\,d}}-1\right)
\end{array} 
\label{NRJ}
\end{equation}
satisfies
$$
\frac{\textstyle d\,E(\boldsymbol{U},\,t)}{\textstyle d\,t}=0,\quad E(\boldsymbol{U},\,t) \geq 0, \quad E(\boldsymbol{U},\,t)=0 \Leftrightarrow \boldsymbol{U}(x,\,t)=\boldsymbol{0}.
$$
\label{LemmeNRJ}
\end{lemme}

{\sc Proof.} We multiply the first equation of (\ref{LCscal}) by $v$, and we integrate it by parts on $\Omega_0$; then the second equation of (\ref{LCscal}) gives
$$
\begin{array}{lll}
\displaystyle
\int_{-\infty}^{\alpha^-}\rho\,v\,\frac{\textstyle \partial\,v}{\textstyle \partial\,t}\,dx&=&
\displaystyle \int_{-\infty}^{\alpha^-}v\,\frac{\textstyle \partial\,\sigma}{\textstyle \partial\,x}\,dx,\\
[10pt]
&=& \displaystyle v(\alpha^-,\,t)\,\sigma(\alpha^-,\,t)-\int_{-\infty}^{\alpha^-}\sigma\,\frac{\textstyle \partial\,v}{\textstyle \partial\,x}\,dx,\\
[10pt]
&=& \displaystyle v(\alpha^-,\,t)\,\sigma(\alpha^-,\,t)-\int_{-\infty}^{\alpha^-}\frac{\textstyle 1}{\textstyle \rho\,c^2}\,\sigma\,\frac{\textstyle \partial\,\sigma}{\textstyle \partial\,t}\,dx.
\end{array}
$$
In the same way, we obtain
$$
\int_{\alpha^+}^{+\infty}\rho\,v\,\frac{\textstyle \partial\,v}{\textstyle \partial\,t}\,dx =
-v(\alpha^+,\,t)\,\sigma(\alpha^+,\,t)-\int_{\alpha^+}^{+\infty}\frac{\textstyle 1}{\textstyle \rho\,c^2}\,\sigma\,\frac{\textstyle \partial\,\sigma}{\textstyle \partial\,t}\,dx.
$$
Adding the two previous equations gives: $\varepsilon_1+\varepsilon_2=0$, where
$$
\begin{array}{lll}
\varepsilon_1 &=& \displaystyle \int_{-\infty}^{+\infty}\rho\,v\,\frac{\textstyle \partial\,v}{\partial\,t}\,dx+\int_{-\infty}^{+\infty}\frac{\textstyle 1}{\textstyle \rho\,c^2}\,\sigma\,\frac{\textstyle \partial\,\sigma}{\textstyle \partial\,t}\,dx\\
[10pt]
&=& \displaystyle \frac{\textstyle d}{\textstyle d\,t}\,\underbrace{
\frac{\textstyle 1}{\textstyle 2}\int_{-\infty}^{+\infty}\left(\rho\,v^2+\frac{\textstyle 1}{\textstyle \rho\,c^2}\,\sigma^2\right)dx}_{\displaystyle E_1},
\end{array}
$$
and, using (\ref{JCUt}),
$$
\begin{array}{lll}
\varepsilon_2 &=& \displaystyle \sigma(\alpha^-,\,t)\,\left[v(\alpha,\,t)\right],\\
[10pt]
&=& \displaystyle \frac{\textstyle 1}{\textstyle K}\,\frac{\textstyle \sigma(\alpha^-,\,t)}{\textstyle \left(\displaystyle 1-\frac{\textstyle \sigma(\alpha^-,\,t)}{\textstyle K\,d}\right)^2}\,\frac{\textstyle \partial\,\sigma}{\textstyle \partial \,t}(\alpha^-,\,t),\\
[10pt]
&=& \displaystyle \frac{\textstyle d}{\textstyle d\,t}\,\underbrace{K\,d^2\left(\ln\left(1-\frac{\textstyle \sigma(\alpha^-,\,t)}{\textstyle K\,d}\right)+\frac{\textstyle 1}{\textstyle 1 - \displaystyle \frac{\textstyle \sigma(\alpha^-,\,t)}{\textstyle K\,d}}-1\right)}_{\displaystyle E_2}.
\end{array}
$$
The quantity $E=E_1+E_2$ therefore satisfies $\frac{d\,E}{d\,t}=0$; $E_1$ is obviously a positive definite quadratic form. All we have to do now is to study the sign of $E_2$. Setting $\theta=1-\sigma(\alpha^-,\,t)/(K\,d)$, a study of $g(\theta)=\ln \theta + 1/\theta-1$ for $\theta >0$ shows that $g(\theta)\geq 0$, and that $g(\theta)=0 \Leftrightarrow \theta=1$, i.e. $\sigma(\alpha^-,\,t)=0 $. From (\ref{JCBB}) and (\ref{SUpm}), one sees that $\sigma(\alpha^-,\,t)=0 \Leftrightarrow \sigma(\alpha^+,\,t)$ and $v(\alpha^\pm,\,t)=0$. $\qquad \Box$

Lemma \ref{LemmeNRJ} means that (\ref{NRJ}) is an energy which is split into two terms. The first term, $E_1$, is the classical mechanical energy associated with the propagation of elastic waves outside the fracture. The second term, $E_2$, is the mechanical energy associated with the nonlinear deformation of the fracture. Since the inertial effects are neglected in (\ref{JCBB}), $E_2$ amounts to a potential energy. Note that in the limit case $|\sigma(\alpha^-,\,t)|\ll K\,d$, one gets $E_2\rightarrow\frac{1}{2\,K}\,\sigma^2(\alpha^-,\,t)$, which corresponds to the well-known potential energy of a linear spring.

\subsection{Method of characteristics}

To express $\boldsymbol{U}(x,\,t)$ in terms of limit-values of the fields at $\alpha$, we now use the {\it Riemann invariants} $J^{R,L}$ that are constant along the {\it characteristics} $\gamma_{R,L}$ \cite{RAVIART96}. The invariants for linear PDE's with constant coefficients are very simple
\begin{equation}
\left\{
\begin{array}{l}
\displaystyle \gamma_R:\,\frac{\textstyle d\,x}{\textstyle d\,t}=+c\,\Rightarrow\,\left.\frac{\textstyle dJ^R}{\textstyle d\,t}\right|_{\gamma_R}=0,\quad \mbox{ with }J^R(x,\,t)=\frac{\textstyle 1}{\textstyle 2}\left(v-\frac{\textstyle 1}{\textstyle \rho\,c}\sigma\right)(x,\,t),\\
[10pt]
\displaystyle \gamma_L:\,\frac{\textstyle d\,x}{\textstyle d\,t}=-c\,\Rightarrow\,\left.\frac{\textstyle dJ^L}{\textstyle d\,t}\right|_{\gamma_L}=0,\quad \mbox{ with }J^L(x,\,t)=\frac{\textstyle 1}{\textstyle 2}\left(v+\frac{\textstyle 1}{\textstyle \rho\,c}\sigma\right)(x,\,t).
\end{array}
\right.
\label{Jrl}
\end{equation}
After some calculations, the first equation of (\ref{JCBB}) and (\ref{Jrl}) along with the initial data condition of compact support in $\Omega_0$ give for $t\geq t_0$
\begin{equation}
\begin{array}{l}
\sigma(\alpha^\pm,\,t)=-\rho_1\,c_1\,v(\alpha^+,\,t),\\
[8pt]
\displaystyle
v(\alpha^-,\,t)=-\frac{\textstyle \rho_1\,c_1}{\textstyle \rho_0\,c_0}\,v(\alpha^+,\,t)+2\,J_0^R\left(\alpha-c_0(t-t_0),\,t_0\right),
\end{array}
\label{SUpm} 
\end{equation}
where the subscript $i$ on $J^{R,L}_i$ refers to $\Omega_i$. The following lemma expresses $\boldsymbol{U}(x,\,t)$ in terms of $v(\alpha^+,\,s)$, with $t_0\leq s \leq t$, and of the initial values of the Riemann invariants, which are linked to the initial data $\boldsymbol{U}_0(x)$ via (\ref{MatA}), (\ref{IBVP}) and (\ref{Jrl}). 
\begin{lemme}
Setting
$$
t_A =t-\frac{\textstyle 1}{\textstyle c_0}(\alpha-x), \qquad
t_B=t-\frac{\textstyle 1}{\textstyle c_1}(x-\alpha),
$$
the solution $\boldsymbol{U}(x,\,t)$ of (\ref{IBVP}) is given by
\begin{equation}
\begin{array}{l}
\underline{x<\alpha}: 
\boldsymbol{U}(x,\,t)=
\left(
\begin{array}{cc}
1           &  1\\
[4pt]
-\rho_0\,c_0 & \rho_0\,c_0
\end{array}
\right)\,
\left(
\begin{array}{c}
J_0^R\left(x-c_0(t-t_0),\,t_0\right)\\
[4pt]
\Delta_A(x,\,t)
\end{array}
\right),\\
\\
\displaystyle \mbox{with } \Delta_A(x,\,t)=
\left\{
\begin{array}{l}
\displaystyle
-\frac{\textstyle \rho_1\,c_1}{\textstyle \rho_0\,c_0}\,v(\alpha^+,\,t_A)+J_0^R(\alpha-c_0(t_A-t_0),\,t_0)
\,\mbox{ if }\,t_A \geq t_0,\\
[8pt]
J_0^L(x+c_0(t-t_0),\,t_0) \, \mbox{ otherwise},
\end{array}
\right.
\\ 
\\
\underline{x>\alpha}: 
\boldsymbol{U}(x,\,t)=
\left(
\begin{array}{c}
1\\
[4pt]
-\rho_1\,c_1
\end{array}
\right)
\Delta_B(x,\,t),\\
\\
\displaystyle \mbox{with } \Delta_B(x,\,t)=
\left\{
\begin{array}{l}
\displaystyle
v(\alpha^+,\,t_B)\, \mbox{ if }\,t_B \geq t_0,\\
[4pt]
0 \, \mbox{ otherwise}.
\end{array}
\right.
\label{U(x,t)}
\end{array}%
\end{equation}
\label{LemmeU(x,t)}
\end{lemme}
{\sc Proof}.
The properties of Riemann invariants (\ref{Jrl}) in $\Omega_0$ lead to
$$
\begin{array}{l}
J_0^R(x,\,t)=J_0^R(x-c_0(t-t_0),\,t_0),\\
J_0^L(x,\,t)=
\left\{
\begin{array}{l}
\displaystyle
J_0^L(\alpha^-,\,t_A) \,\mbox{ if } t_A \geq t_0,\\
J_0^L(x+c_0(t-t_0),\,t_0) \,\mbox{ otherwise}.
\end{array}
\right.
\end{array}
$$
For $(x,\,t)$ such that $t_A \geq t_0$, (\ref{Jrl}) and (\ref{SUpm}) imply that
$$
J_0^L(\alpha^-,\,t_A)=-\frac{\textstyle \rho_1\,c_1}{\textstyle \rho_0\,c_0}\,v(\alpha^+,\,t_A)+J_0^R(\alpha-c_0(t_A-t_0),\,t_0).
$$
The value of $\Delta_A$ and the previous equations allow to conclude for $x<\alpha$. Since the support of $\boldsymbol{U}_0(x)$ is included in $\Omega_0$, (\ref{Jrl}) in $\Omega_1$ lead to
$$
\begin{array}{l}
J_1^R(x,\,t)=
\left\{
\begin{array}{l}
\displaystyle
J_1^R(\alpha^+,\,t_B)\,\mbox{ if } t_B \geq t_0,\\
0\,\mbox{ otherwise},
\end{array}
\right.
\\
J_1^L(x,\,t)=J_1^L(x+c_1(t-t_0),\,t_0)=0.
\end{array}
$$
For $(x,\,t)$ such that $t_B \geq t_0$, (\ref{Jrl}) and (\ref{SUpm}) imply that
$$
J_1^R(\alpha^+,\,t_B)=v(\alpha^+,\,t_B).
$$
The value of $\Delta_B$ and the previous equations allow to conclude for $x>\alpha$. $\Box$
\begin{lemme}
The limit value $y=v(\alpha^+,\,t)$ satisfies the nonlinear ODE 
\begin{equation}
\left|
\begin{array}{l}
\displaystyle
\frac{\textstyle d\,y}{\textstyle d\,t}=f(y,\,t),\qquad y(t_0)=0,\quad \mbox{ with}\\
[8pt]
\displaystyle
f(y,\,t)=\frac{\textstyle K}{\textstyle \rho_1\,c_1}\left(1+\frac{\textstyle \rho_1\,c_1}{\textstyle K\,d}\,y\right)^2\left(g(t)-\left(1+\frac{\textstyle \rho_1\,c_1}{\textstyle \rho_0\,c_0}\right)y\right),\\
[8pt]
\displaystyle
g(t)=2\,J_0^R\left(\alpha-c_0\,(t-t_0),\,t_0\right).
\end{array}
\right.
\label{ODE}
\end{equation}
\label{LemmeODE}
\end{lemme}
{\sc Proof}. From (\ref{JCBBbis}) and the first equation of (\ref{SUpm}), we deduce 
$$
\begin{array}{lll}
v(\alpha^+,\,t)-v(\alpha^-,\,t) &=& \displaystyle \frac{\textstyle 1}{\textstyle K}\,\frac{\textstyle 1}{\textstyle \left(\displaystyle 1-\frac{\textstyle \sigma(\alpha^+,\,t}{\textstyle K\,d}\right)^2}\,\frac{\textstyle \partial\,\sigma}{\textstyle \partial \,t}(\alpha^+,\,t),\\
&=& -\displaystyle \frac{\textstyle \rho_1\,c_1}{\textstyle K}\,\frac{\textstyle 1}{\textstyle \left(\displaystyle 1+\frac{\textstyle \rho_1\,c_1}{\textstyle K\,d}\,v(\alpha^+,\,t)\right)^2}\,\frac{\textstyle \partial\,v}{\textstyle \partial \,t}(\alpha^+,\,t).
\end{array}
$$
$v(\alpha^-,\,t)$ is then eliminated via the second equation of (\ref{SUpm}), giving (\ref{ODE}). $\qquad \Box$

\subsection{Well-posedness of the initial-boundary value problem}

As shown by lemmas \ref{LemmeU(x,t)} and \ref{LemmeODE}, the solution of (\ref{IBVP}) is expressed in a unique manner in terms of $v(\alpha^+,\,t)$ solution of the ODE (\ref{ODE}). Proving the existence and uniqueness of the solution to (\ref{IBVP}) therefore amounts to showing that the solution to (\ref{ODE}) exists and is unique, as followed. 
\begin{theorem}
Let $\boldsymbol{U}_0(x) \in C^p_c(\mathbb{R})$, $p\geq 1$. There exists a unique global solution $\boldsymbol{U}(x,\,t)\in C^p\left(\mathbb{R}\setminus \alpha,\,[t_0,\,+\infty[\right)$ to the initial-boundary value problem (\ref{IBVP}). 
\label{ThExistence}
\end{theorem}
{\sc Proof}.
For $t>t_0$, $y \rightarrow f(y,\,t)$ in (\ref{ODE}) is $C^\infty$, it is therefore a locally Lipschitz function. Moreover, $J^{R}_0(x,\,t_0)$ is a $C^p_c$ function ($p\geq 1$) deduced from $\boldsymbol{U}_0(x)$, therefore $g(t)$ in (\ref{ODE}) is $C^p_c$: hence $t \rightarrow f(y,\,t)$ is continuous. The Cauchy-Lipschitz theorem ensures that the solution $y(t)$ is unique, if it exists. Lastly, $f\in C^p$ implies that $y\in C^{p+1}$ \cite{DEMAILLY}: the lemma \ref{LemmeU(x,t)} ensures the $C^p$ regularity of $\boldsymbol{U}$.

Suppose that $y(t)$ is not bounded as $t \rightarrow t^*$; the first equation of (\ref{SUpm}) implies that $\sigma(\alpha^-,\,t)\rightarrow \pm\infty$. Since $\sigma(\alpha^-,\,t)<K\,d$ (see section \ref{SecJC0}), only the case $\sigma(\alpha^-,\,t)\rightarrow -\infty$ needs to be addressed. In this case, (\ref{NRJ}) implies that $E(\boldsymbol{U},\,t)\rightarrow + \infty$, which is impossible: lemma \ref{LemmeNRJ} and $\boldsymbol{U}_0 \in C^p_c$ mean that $E(\boldsymbol{U},\,t)=E(\boldsymbol{U}_0,\,t_0)<+\infty$. Hence $y(t)$ is always bounded, and the local existence due to the Cauchy-P\'eano theorem is also global \cite{CROUZEIX}. $\qquad \Box$

\section{Numerical treatment}\label{SecNum}

\subsection{Numerical scheme}\label{SecSchema}

\begin{figure}[htbp]
\begin{center}
\includegraphics[scale=0.9]{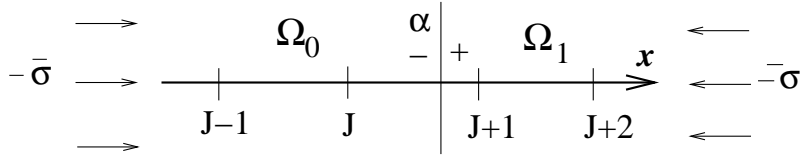}
\end{center}
\caption{1D rock fractured at $x=\alpha$; spatial mesh.}
\label{FigTampon1D}
\end{figure}

Given $(x_i,\,t_n)=(i\,\Delta\,x,\,t_0 + n\,\Delta\,t)$, where $\Delta\,x$ is the mesh size and $\Delta \,t$ is the time step, we seek an approximation $\boldsymbol{U}_i^n$ of $\boldsymbol{U}(x_i,\,t_n)$. We use two-step, explicit, and $(2\,s+1)$-point spatially-centered finite-difference schemes to integrate (\ref{IBVP}). Time-stepping can then be written symbolically 
\begin{equation}
\boldsymbol{U}_i^{n+1} = \boldsymbol{U}_i^n+\boldsymbol{H}_q\left(\boldsymbol{U}_{i-s}^n,...,\,\boldsymbol{U}_{i+s}^n\right),
\label{TM}
\end{equation}
with $q=0$ if $x_i\in\Omega_0$, $q=1$ otherwise, and with $\boldsymbol{H}_q:\mathbb{R}^{2 \times (2s+1)}\rightarrow \mathbb{R}^2$ \cite{RAVIART96}. We define $J$ so that $x_J \leq \alpha < x_{J+1}$ (figure \ref{FigTampon1D}). A grid point is {\it regular} if all the grid points $x_{i-s},\,...,\,x_{i+s}$ used in (\ref{TM}) belong to the same medium as $x_i$; in this case, (\ref{TM}) is applied classically. Otherwise, a grid point is {\it irregular}, and its time-stepping is described in the next section. The irregular points are $x_{J-s+1},...\,,x_{J+s}$. 

For the numerical experiments in section \ref{SecExp}, we choose the ADER $r$ schemes \cite{ADER04}, where $r$ denotes the order of accuracy. The stability limit is CFL = $\max(c)\,\Delta\,t\,/\,\Delta\,x=1$. For odd values of $r$, as chosen here, $s=r\,/\,2$. Other schemes can be used (numerical experiments have been performed successfully with a  flux-limiter scheme and a fifth-order WENO scheme): a priori, readers can adapt their favorite solver to the forthcoming discussion.

\subsection{Interface method}\label{SecESIM}
 
One applies the explicit simplified interface method (ESIM) \cite{ALIMENTAIRE1,BIBLE1}: at irregular points, some of the numerical values used for the time-stepping procedure (\ref{TM}) are modified. At time $t_n$ and at irregular points $x_i$, these {\it modified values} $\boldsymbol{U}_i^*$ are numerical estimates of smooth extension $\boldsymbol{U}^*(x,\,t_n)$ of $\boldsymbol{U}(x,\,t_n)$ across $\alpha$ 
\begin{equation}
\begin{array}{l}
\displaystyle
x>\alpha,\quad  \boldsymbol{U}^{*}(x,\,t_n) = \sum_{m=0}^{2\,k-1}\frac{\textstyle(x-\alpha)^m}{\textstyle m\,!} \frac{\textstyle \partial^m}{\textstyle \partial\, x^m}\,\boldsymbol{U}(\alpha^-,\,t_n),\\
[15pt]
\displaystyle
x\leq \alpha,\quad  \boldsymbol{U}^{*}(x,\,t_n) = \sum_{m=0}^{2\,k-1}\frac{\textstyle( x-\alpha)^m}{\textstyle m\,!} \frac{\textstyle \partial^m}{\textstyle \partial\, x^m}\,\boldsymbol{U}(\alpha^+,\,t_n).
\end{array}
\label{Umod}
\end{equation}
The next paragraph details how to estimate $\frac{\partial^m}{\partial\, x^m}\,\boldsymbol{U}(\alpha^\pm,\,t_n)$ in (\ref{Umod}). 

{\bf Estimation of $\frac{\partial^m}{\partial\, x^m}\,\boldsymbol{U}(\alpha^-,\,t_n)$}. We focus here on the computation of the modified values at irregular points $x_j$ on $\Omega_1$ ($j=J+1,...,J+s$); the case of irregular points on $\Omega_0$ ($j=J-s+1,...,J$) is similar. To estimate $\frac{\partial^m}{\partial\, x^m}\,\boldsymbol{U}(\alpha^-,\,t_n)$ in (\ref{Umod}), we write Taylor expansions of $\boldsymbol{U}(x_i,\,t_n)$ on the left of $\alpha$ ($i=J-k+1,\,...,\,J$)
\begin{equation}
\boldsymbol{U}(x_i,t_n) =\sum_{m=0}^{2\,k-1} \frac{\textstyle (x_i -\alpha)^m}{\textstyle m\,!}\, \frac{\textstyle \partial^m}{\textstyle \partial \,x^m}\,\boldsymbol{U}(\alpha^-,\,t_n)+\boldsymbol{O}(\Delta\, x ^{2k}),
\label{TaylorM}
\end{equation}
and on the right of $\alpha$ ($i=J+1,\,...,\,J+k$), with the jump conditions (\ref{JCBBm})
\begin{equation}
\begin{array}{l}
\boldsymbol{U}(x_i,t_n) = \displaystyle \sum_{m=0}^{2\,k-1} \frac{\textstyle (x_i -\alpha)^m}{\textstyle m\,!}\frac{\textstyle \partial^m}{\textstyle \partial x^m}\boldsymbol{U}(\alpha^+,t_n)+\boldsymbol{O}(\Delta \,x ^{2k})\\
[15pt]
\quad=\displaystyle \sum_{m=0}^{2\,k-1} \frac{\textstyle (x_i -\alpha)^m}{\textstyle m\,!}\boldsymbol{D}_m\left(\boldsymbol{U}(\alpha^-,t_n) ,...,
\frac{\textstyle \partial^{m}}{\textstyle \partial\,x^{m}}\,\boldsymbol{U}(\alpha^-,t_n),\, 
\frac{\textstyle \partial^{m+1}}{\textstyle \partial\,x^{m+1}}\,\boldsymbol{U}(\alpha^-,t_n) \right)\\
\quad +\boldsymbol{O}(\Delta \,x ^{2k}).
\end{array}
\label{TaylorP}
\end{equation}
The exact values $\boldsymbol{U}(x_i,\,t_n)$ are replaced by the known $\boldsymbol{U}_i^n$ in (\ref{TaylorM}) and (\ref{TaylorP}); $\frac{\partial^{2\,k}}{\partial \,x^{2\,k}}\,\boldsymbol{U}(\alpha^-,t_n)$ and the Taylor remainders are also removed, leading to 
\begin{equation}
\boldsymbol{F}\left(\boldsymbol{U}^-,...,\frac{\textstyle \partial^{2\,k-1}}{\textstyle \partial\,x^{2\,k-1}}\,\boldsymbol{U}^-\right)=\left(\boldsymbol{U}_{J-k+1}^n,\,...,\,\boldsymbol{U}_{J+k}^n \right)^T,
\label{SysNl}
\end{equation}
where $\frac{\partial^{m}}{\partial \,x^{m}}\,\boldsymbol{U}^-$ are estimates of $\frac{\partial^{m}}{\partial \,x^{m}}\,\boldsymbol{U}(\alpha^-,t_n)$, and where $\boldsymbol{F}:\mathbb{R}^{4\,k}\rightarrow\mathbb{R}^{4\,k}$ is a nonlinear application. The system (\ref{SysNl}) is solved using  Newton's method. 

{\bf Computation and use of modified values}.
Once the limit values $\frac{\partial^{m}}{\partial\,x^{m}}\,\boldsymbol{U}^-$ are known, the modified values $\boldsymbol{U}_i^*$ are deduced from $\boldsymbol{U}^*(x_i,\,t_n)$ (\ref{Umod}):
\begin{equation}
i=J+1,...,J+s,\quad \boldsymbol{U}^{*}_i = \sum_{m=0}^{2\,k-1}\frac{\textstyle(x_i-\alpha)^m}{\textstyle m\,!} \frac{\textstyle \partial^m}{\textstyle \partial\, x^m}\,\boldsymbol{U}^-.
\label{U*}
\end{equation}
When all the modified values have been computed, all that remains to be done is to perform time-stepping at irregular points. To do so, the time-stepping (\ref{TM}) at an irregular point $x_i$ requires the use of the numerical values at the grid points in the same medium as $x_i$, and the modified values otherwise
\begin{equation}
\begin{array}{l}
i=J-s+1,...,J,\quad \boldsymbol{\boldsymbol{U}}_i^{n+1}=\boldsymbol{\boldsymbol{U}}_i^n+\boldsymbol{H}_0\left(\boldsymbol{\boldsymbol{U}}_{i-s}^n,...,\,\boldsymbol{\boldsymbol{U}}_J^n,\,\boldsymbol{\boldsymbol{U}}_{J+1}^*,...,\,\boldsymbol{\boldsymbol{U}}_{i+s}^*\right),\\
[8pt]
i=J+1,...,J+s,\quad
\boldsymbol{\boldsymbol{U}}_i^{n+1}=\boldsymbol{\boldsymbol{U}}_i^n+\boldsymbol{H}_1\left(\boldsymbol{\boldsymbol{U}}_{i-s}^*,...,\,\boldsymbol{\boldsymbol{U}}_J^*,\,\boldsymbol{\boldsymbol{U}}_{J+1}^n,...,\,\boldsymbol{\boldsymbol{U}}_{i+s}^n\right).
\end{array}
\label{TMESIM}
\end{equation}
{\bf Comments}. We do not propose a numerical analysis of the interface method in the nonlinear context. We just recall some properties that are true in the limit case of negligible nonlinearities; readers are referred to \cite{ALIMENTAIRE1} for proofs. For $|\sigma(\alpha^-,\,t)| / (K\,d)= 0$, the following properies are satisfied:
\begin{description}
\item[{\it (i)}] the system (\ref{SysNl}) has always a unique solution;
\item[{\it (ii)}] in the limit case of a homogeneous medium without any fracture (i.e., $\rho_0 \rightarrow \rho_1$, $c_0 \rightarrow c_1$, $K \rightarrow +\infty$) and if $k\geq s$, then $\boldsymbol{U}_i^* \rightarrow \boldsymbol{U}_i^n$: the interface method (\ref{TMESIM}) amounts to the classical time-marching (\ref{TM});
\item[{\it (iii)}] for a $r$-th order scheme, the local truncation error of (\ref{TMESIM}) is still $r$-th order accurate if $2\,k-1\geq r$.
\end{description}
The convergence measurements performed in section \ref{SecExp} indicate that the property {\it (iii)} is still satisfied in the nonlinear context, but a rigorous proof remains to be obtained. Up to now, we have not chosen $k$: for that purpose, we follow two criteria. First, the spatial derivatives in (\ref{TaylorM}) and (\ref{TaylorP}) need to be well-defined: theorem \ref{ThExistence} implies $2\,k-1\leq p$. Second, we want the properties {\it (i,ii,iii)} to be  satisfied in the linear limit. It leads to the following inequalities
\begin{equation}
\left\{
\begin{array}{l}
2\,k-1 \leq p,\\
[5pt]
k \geq s,\\
[5pt]
2\,k-1 \geq r.
\end{array}
\right.
\label{Koptim}
\end{equation}
In practice, we use the minimum value of $k$ that satisfies (\ref{Koptim}). 

We have no theoretical results about the stability of the interface method, even in the linear context. We have studied many geometrical configurations, values of the physical parameters, and values of $K$, $d$. With a wide range of parameters, no instabilities are usually observed up to the CFL limit, even after very long integration times. However, instabilities are observed when the physical parameters differ considerably between the two sides of an interface. In practice, this is not penalizing when dealing with realistic media. 

The nonlinear system (\ref{SysNl}) may have more than one solution, and we can not be sure that Newton's algorithm selects the right one. To check this point, we compare the numerical value $v^-$ obtained by solving (\ref{SysNl}) with the semi-analytical value deduced from  (\ref{ODE}) and (\ref{SUpm}). For weak to moderate nonlinear effects, the numerical value and the semi-analytical value are the same (up to the accuracy of the integrations). But when the nonlinear effects are large, convergence towards a wrong solution can occur if the mesh used is too coarse. This is not very surprising: since wave profiles tend to be stiffened (see section \ref{SecExp}), the Taylor expansions in (\ref{TaylorM}) and (\ref{TaylorP}) give rather poor estimates. The strategy followed in section \ref{SecExp} consists in using a finer mesh. A better strategy might consist in solving (\ref{SysNl}) with the constraint $\sigma^-<K\,d$, as required by (\ref{JCBB}).

\section{Numerical experiments}\label{SecExp}

We study a 400-m domain fractured at $\alpha=200.67$ m, with parameters \cite{ZHAO01}
$$
\left\{
\begin{array}{l}
\rho_0=\rho_1=1200\,\mbox{ kg/m}^3,\quad K=1.3.10^9\,\mbox{ kg/s}^{2},\\
[5pt]
c_0=c_1=2800\,\mbox{  m/s}, \quad \quad \, d=6.1\,10^{-4}\,\mbox{m}. 
\end{array}
\right.
$$
Since $\rho$ and $c$ are the same on both sides of $\alpha$, the reflected wave and the distorsion of the transmitted wave are induced only by the mechanical behavior of the fracture. Many numerical experiments have been tested successfully with different physical parameters on both sides of $\alpha$. 

\subsection{Test 1: initial-boundary value problem}\label{SecTest1}

\begin{figure}[htbp] 
\begin{center}
\begin{tabular}{cc}
(a) &  (a) \\
\includegraphics[width=6.8cm,height=6cm]{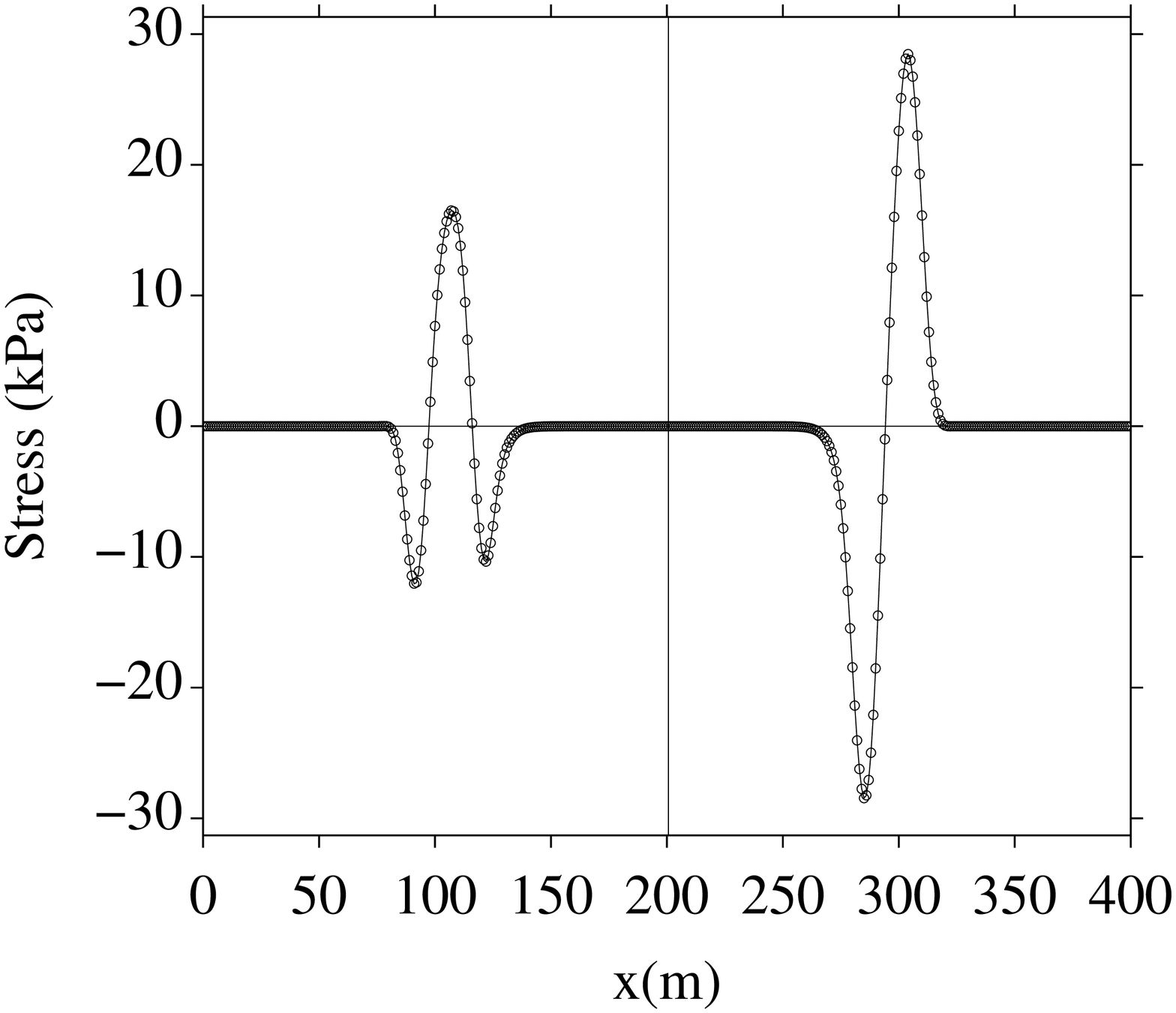}&
\includegraphics[width=6.8cm,height=6cm]{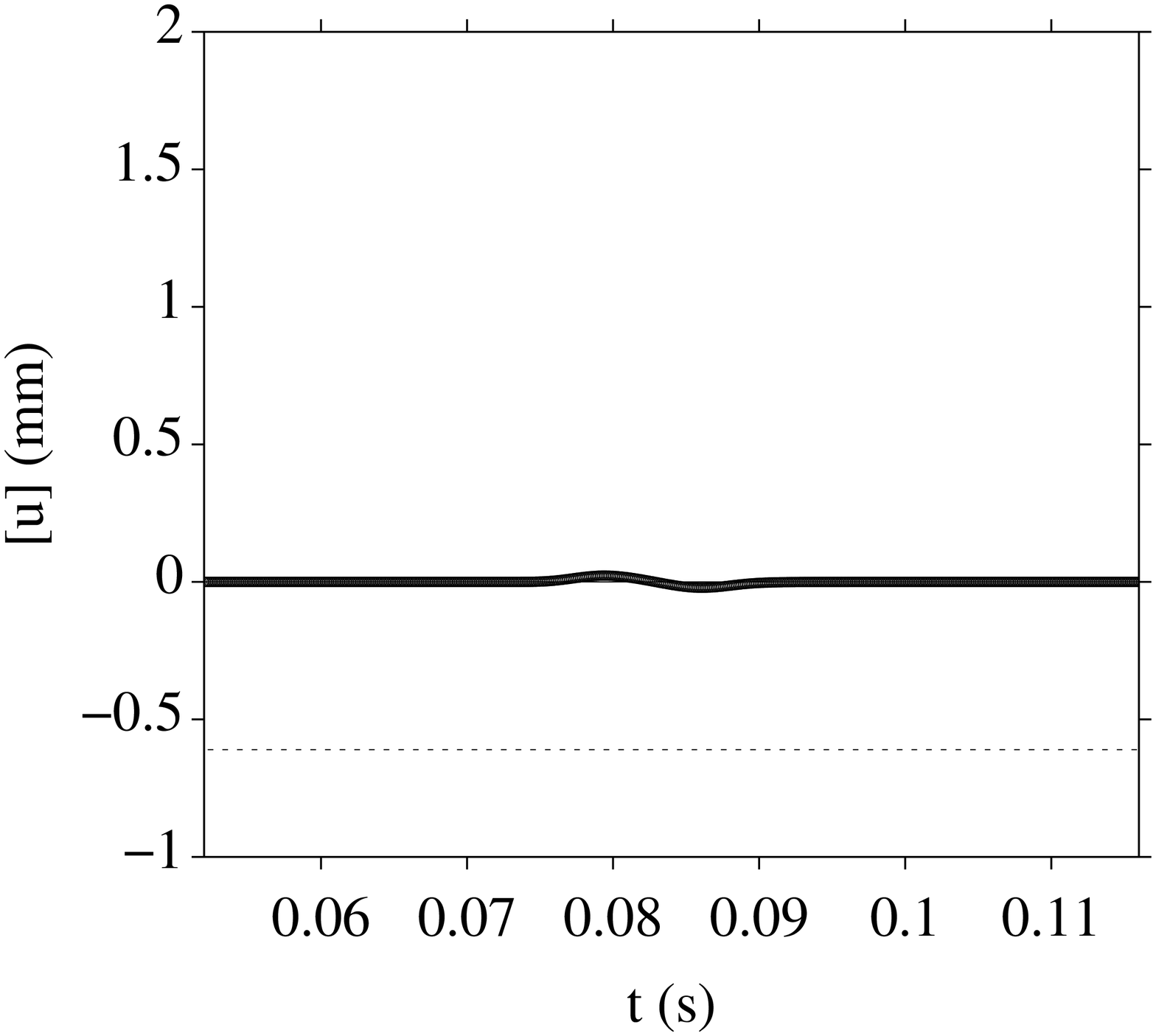}\\
(b)  & (b) \\
\includegraphics[width=6.8cm,height=6cm]{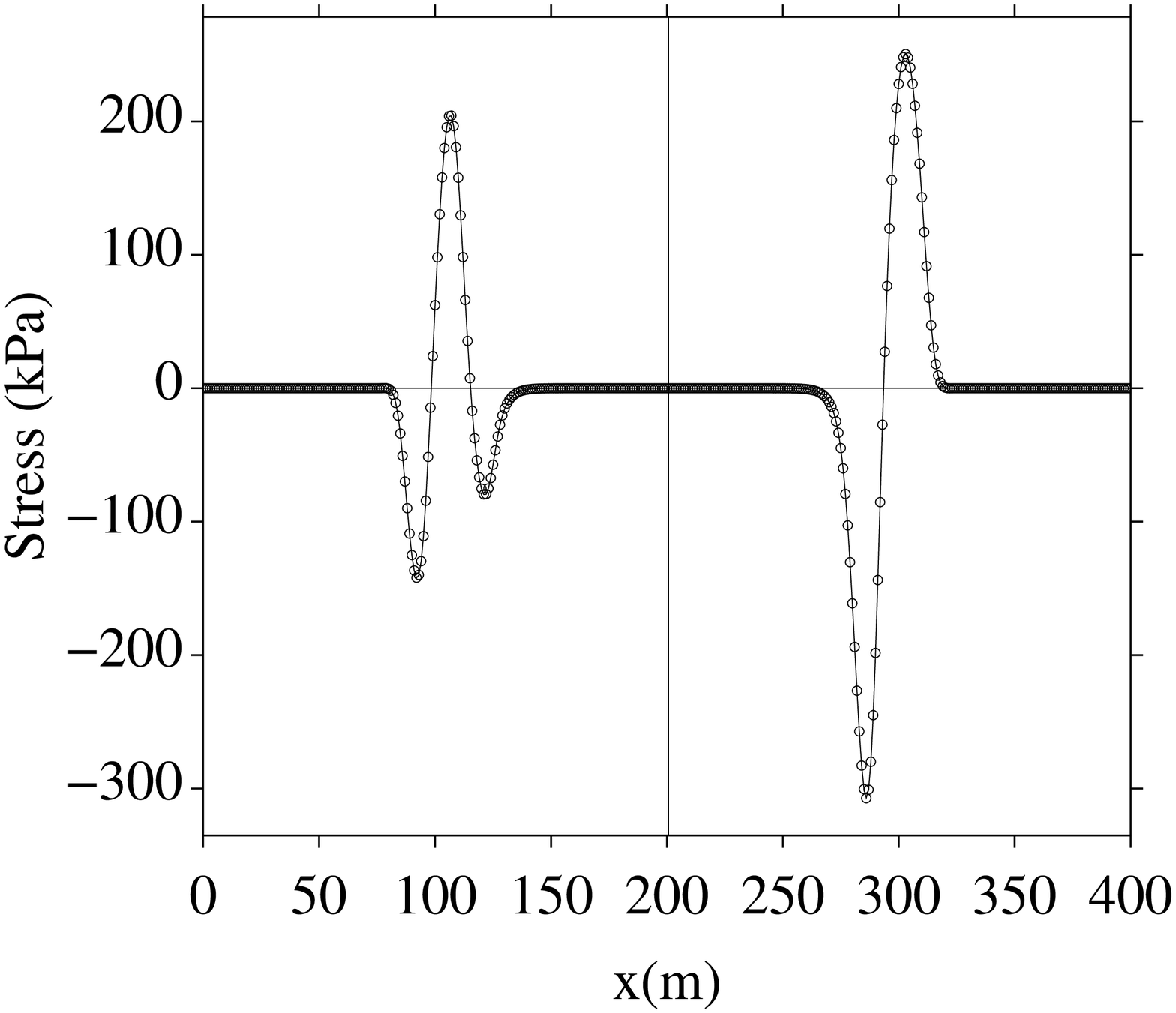}&
\includegraphics[width=6.8cm,height=6cm]{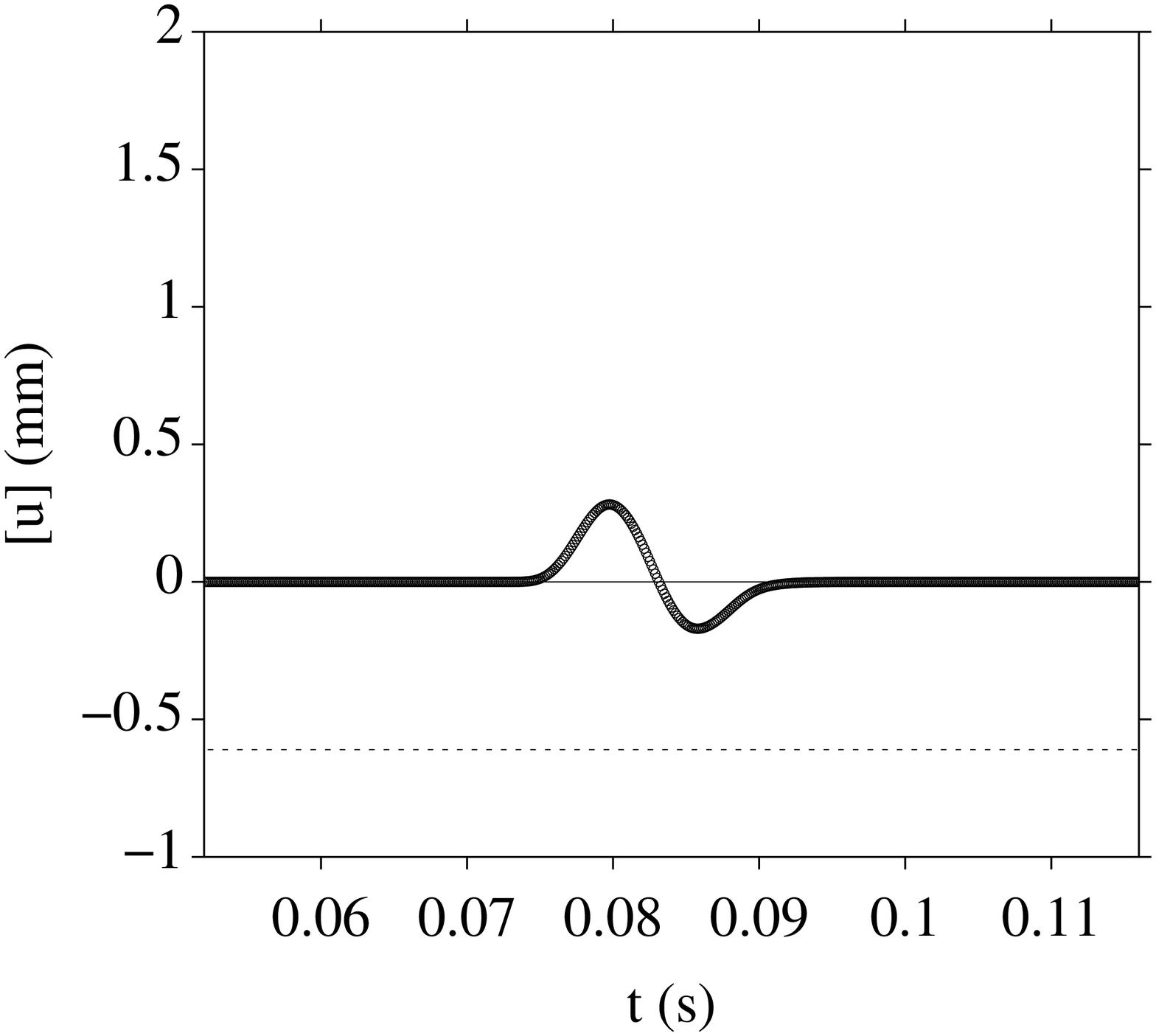}\\
(c) &  (c) \\
\includegraphics[width=6.8cm,height=6cm]{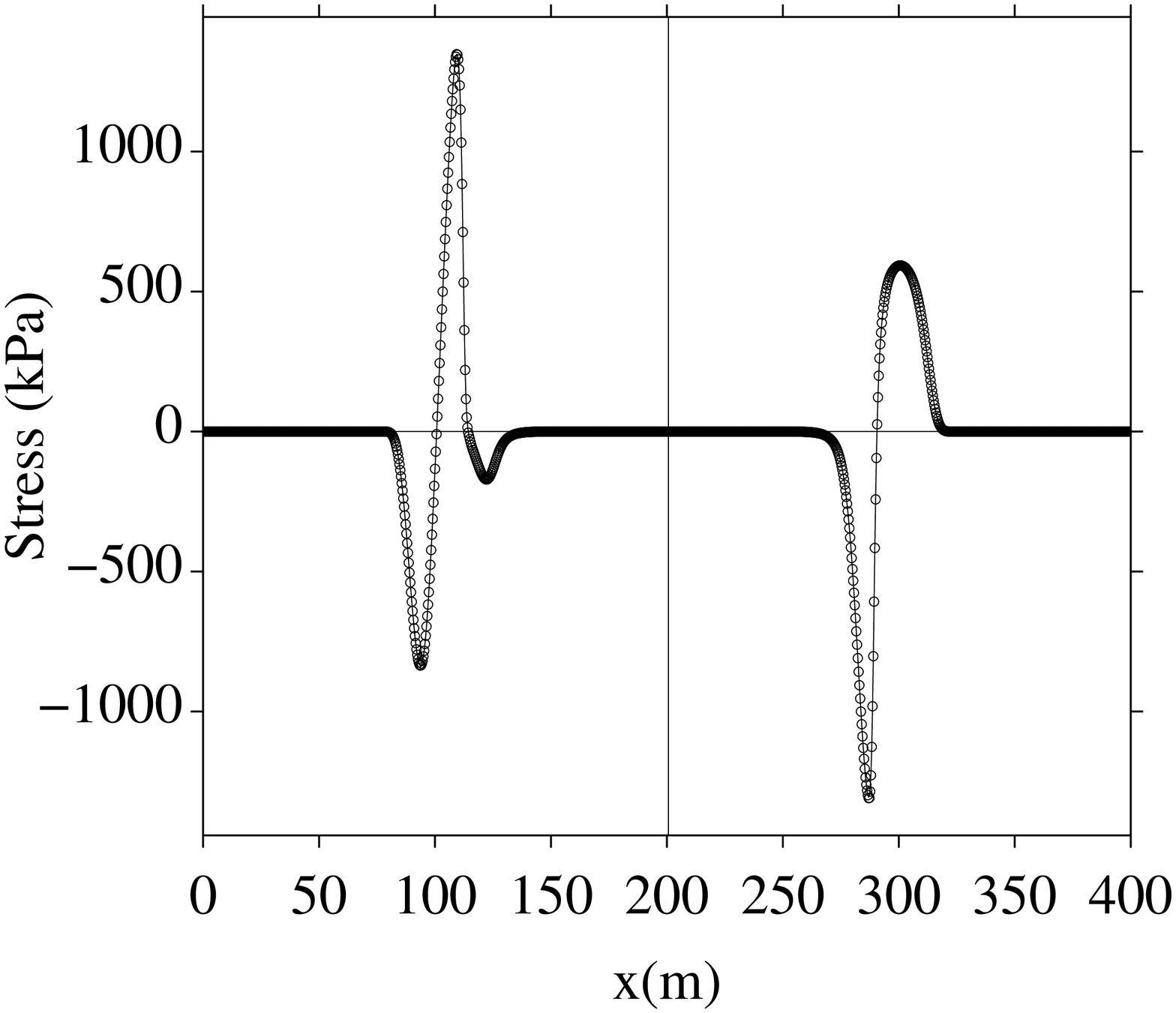}&
\includegraphics[width=6.8cm,height=6cm]{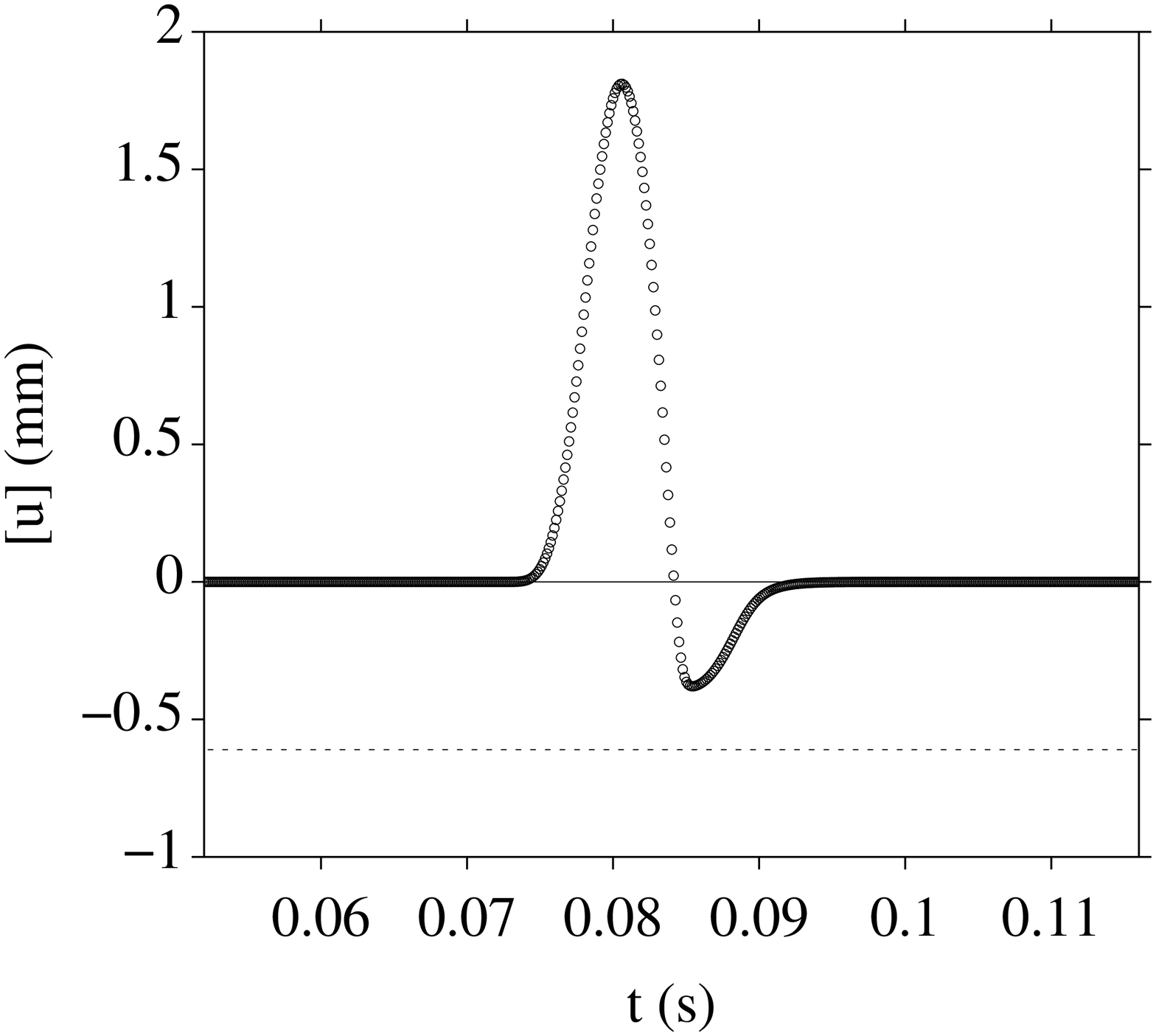}
\end{tabular}
\end{center}
\caption{Test 1: $v_0$ = 0.01 m/s (a), $v_0$ = 0.1 m/s (b), and $v_0$ = 0.4 m/s (c). Left column: numerical (...) and exact (-) values of $\sigma$; right column: numerical values of $[u]$ (the horizontal dotted line denotes $-d$).}
\label{Test1}
\end{figure}

In the first experiment, the initial data is $\boldsymbol{U}_0(x)=(-1/c_0,\,\rho_0)^T\,h(t_0-x\,/\,c_0)$, where $h$ is a spatially-bounded $C^6_c$ combination of sinusoids 
\begin{equation}
h(\xi)=
\left\{
\begin{array}{l}
\displaystyle
\varepsilon\, \displaystyle \sum_{m=1}^4 a_m\,\sin(\beta_m\,\omega_c\,\xi)\quad \mbox{ if  }\, 0<\xi<\frac{\textstyle 1}{\textstyle f_c},\\
[8pt]
0 \,\mbox{ otherwise}, 
\end{array}
\right.
\label{JKPS}
\end{equation}
with $\beta_m=2^{m-1}$, $\omega_c=2\pi\,f_c$; the coefficients $a_m$ are: $a_1=1$, $a_2=-21/32$, $a_3=63/768$, $a_4=-1/512$. The central frequency is $f_c=\omega_c/(2\,\pi)=$ 50 Hz, and $t_0=52$ ms. This choice ensures that the incident wave is a purely rightward-travelling wave, originally located in $\Omega_0$. Three values of $\varepsilon$ are considered, leading to three amplitudes $v_0$ of $v(x,\,t_0)$: 0.01 m/s (a), 0.1 m/s (b), and 0.4 m/s (c). ADER 4 is coupled with the ESIM with $k=3$, as defined by (\ref{Koptim}). The computations are performed with CFL = 0.9 on $N$ = 400 grid points ((a) and (b)) or $N=1200$ grid points (c): here the number of grid points is larger, otherwise (\ref{SysNl}) would not have converged to the right solution.

In the left column of Figure \ref{Test1}, one shows the numerical and analytical values of $\sigma$ at $t$ = 116.29 ms. At this instant, the incident wave has crossed the fracture and one sees the reflected and transmitted waves. In the right column of Figure \ref{Test1}, one shows the time history of the numerical values of $[u(\alpha,\,t)]$ deduced from (\ref{JCBB}) and (\ref{SysNl}). Here the horizontal dotted line represents $-d$, and the vertical scales are the same for (a), (b) and (c). 

In case (a), $v_0$ is too small to mobilize the nonlinearity of the fracture: almost no differences could be detected with simulations performed with (\ref{JClin}). Case (b) corresponds to realistic seismic waves recorded during on-site investigations: moderate nonlinear effects are present. Case (c) corresponds to incident blasting waves: large nonlinear effects are present, stiffening the fronts. In the right column, one clearly sees how $\min([u])$ approaches $-d$ when the incident wave amplitudes increases, without reaching this value as required by (\ref{Penetration}).

\subsection{Test 2: convergence measurements}\label{SecTest2}

\begin{figure}[htbp] 
\begin{center}
\begin{tabular}{cc}
(a) &  (b) \\
\includegraphics[width=6.8cm,height=6cm]{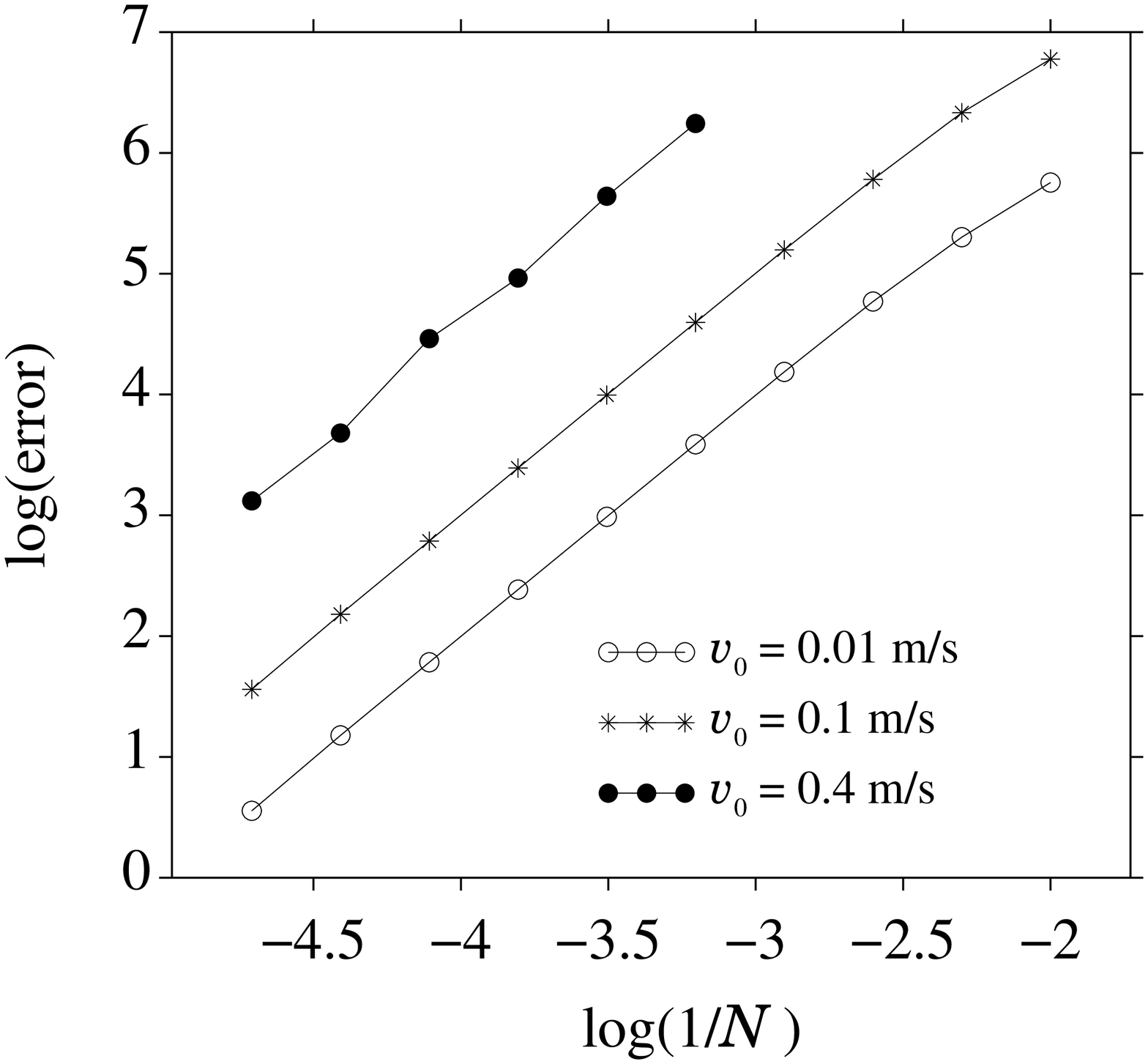}&
\includegraphics[width=6.8cm,height=6cm]{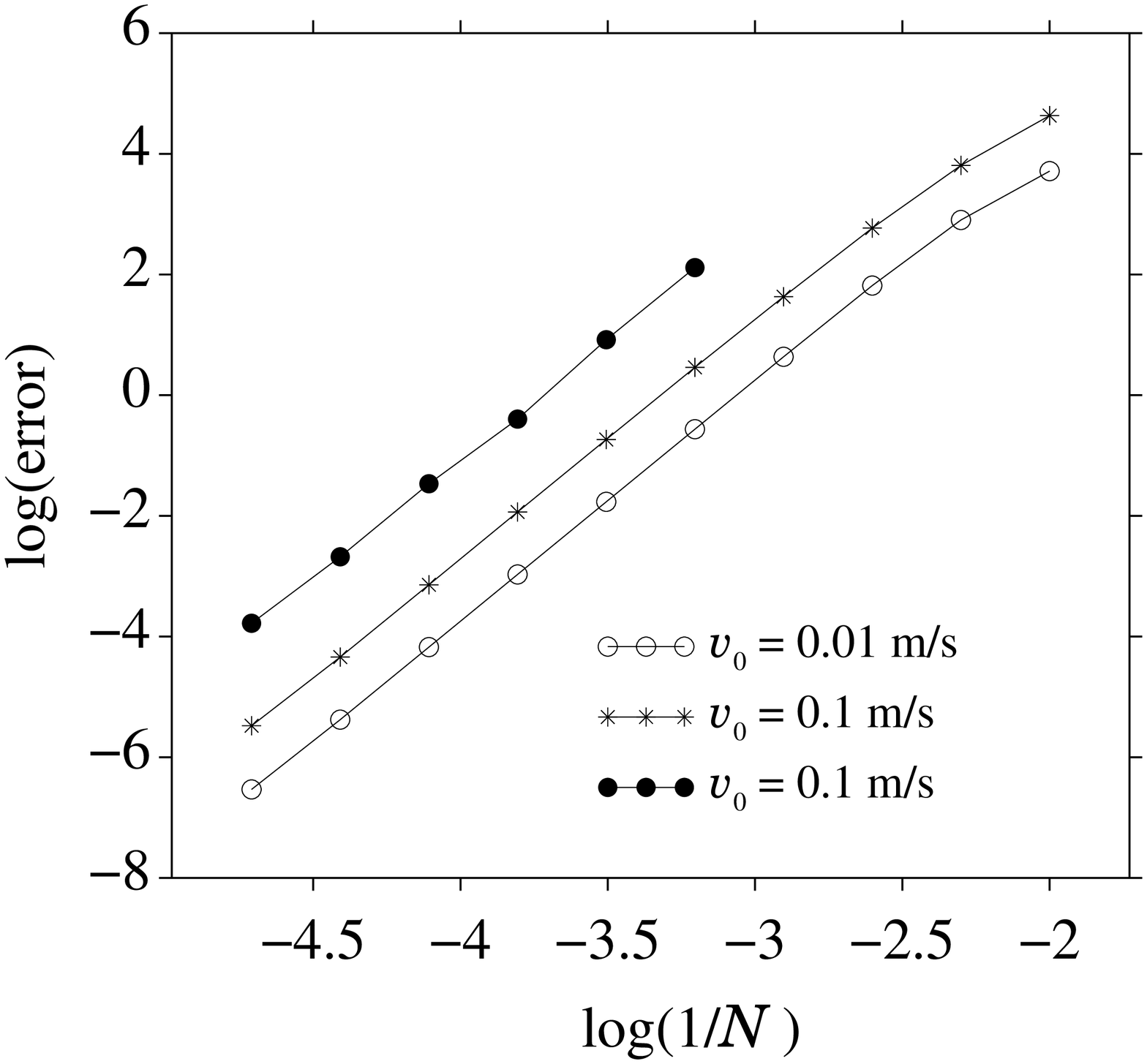}
\end{tabular}
\end{center}
\caption{Test 2: convergence measurements for ADER 2 (a) and ADER 4 (b).}
\label{Test2}
\end{figure}

To check the accuracy of the the interface method, we compare the analytical and numerical values of $\sigma$ on successive refined meshes. The parameters are the same than in subsection \ref{SecTest1}. The errors are measured in norm $L_1$ at $t$ = 116.29 ms. All the computing is carried out in double precision. \newpage

Figure \ref{Test2} shows convergence measurements with ADER 2 (a) and ADER 4 (b) coupled with the ESIM ((\ref{Koptim}) implies $k=2$ for ADER 2). The values of $v_0$ are the same as in test 1: 0.01 m/s (o), 0.1 m/s ($*$), and 0.4 m/s ($\bullet$). The vertical scale of (b) is twice as large than the vertical scale of (a). In all cases, the expected orders of accuracy are exactly reached: order 2 in the case of ADER 2, order 4 in that of ADER 4. These results indicate that coupling a $r$-th order accurate scheme with the interface method is still $r$-th order if the appropriate value of $k$ (\ref{Koptim}) is used.

\subsection{Test 3: sinusoidal source term}\label{Test3}

\begin{figure}[htbp] 
\begin{center}
\begin{tabular}{cc}
(a) &  (a) \\
\includegraphics[width=6.8cm,height=6cm]{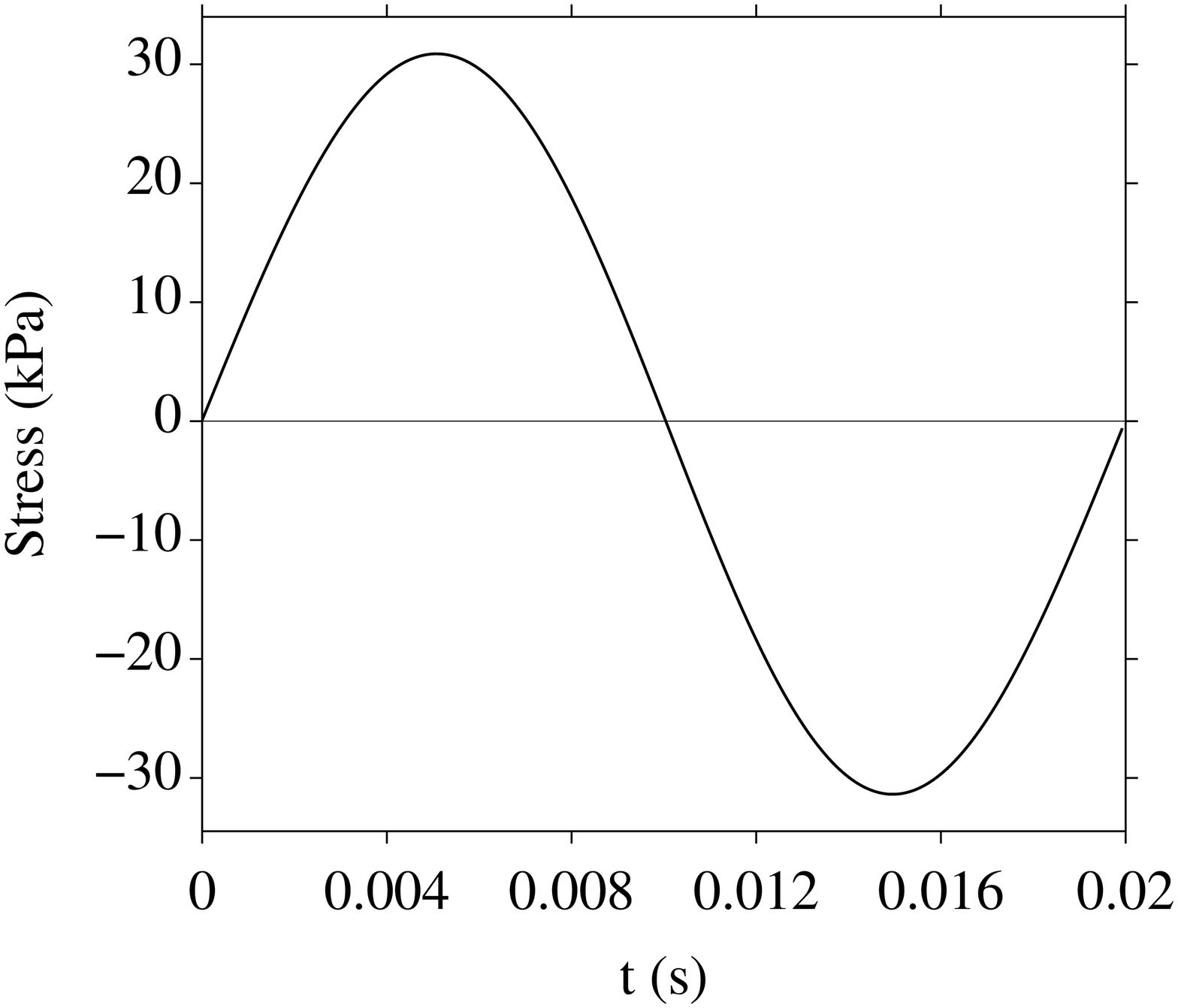}&
\includegraphics[width=6.8cm,height=6cm]{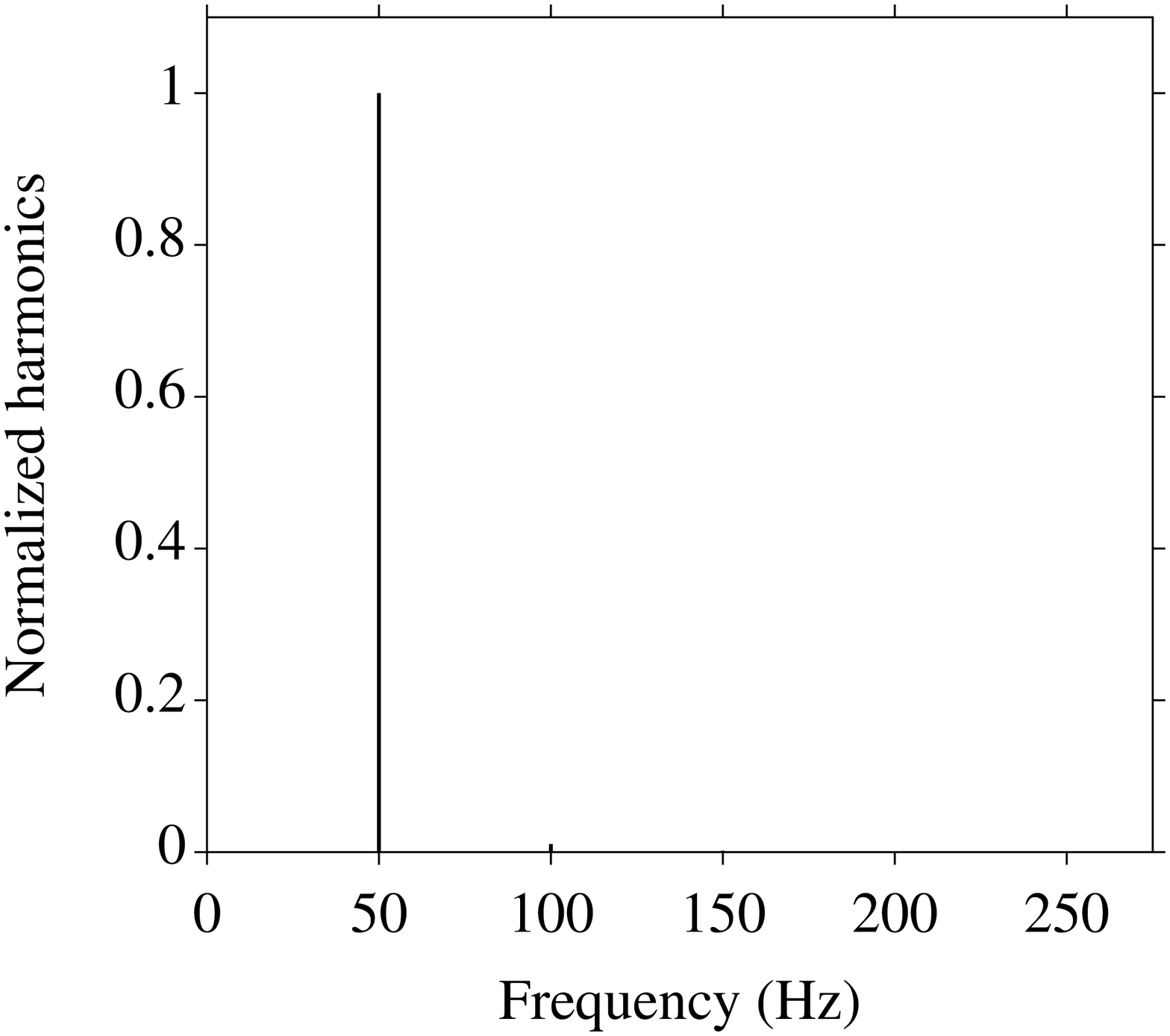}\\
(b) & (b) \\
\includegraphics[width=6.8cm,height=6cm]{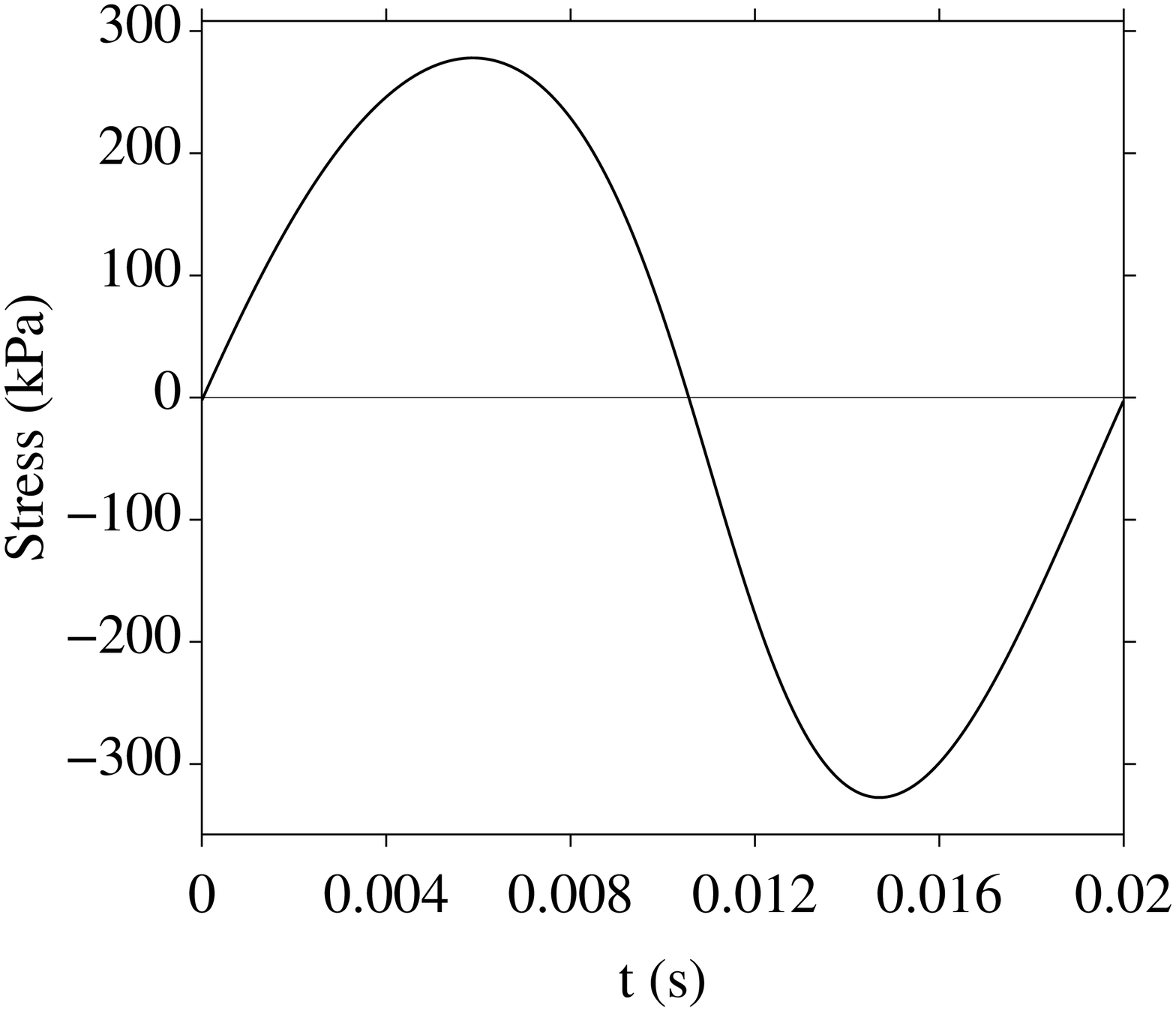}&
\includegraphics[width=6.8cm,height=6cm]{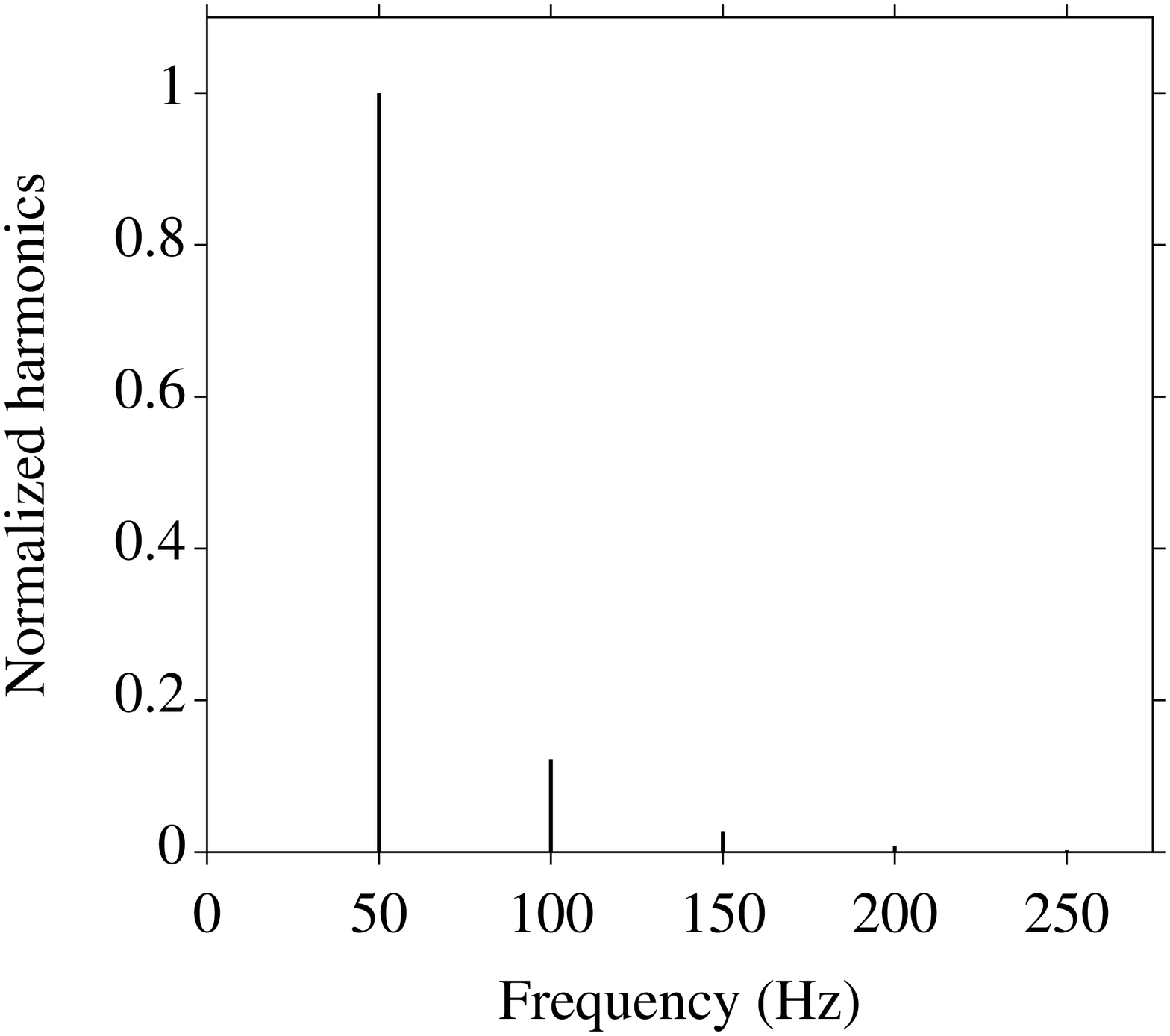}\\
(c) & (c) \\
\includegraphics[width=6.8cm,height=6cm]{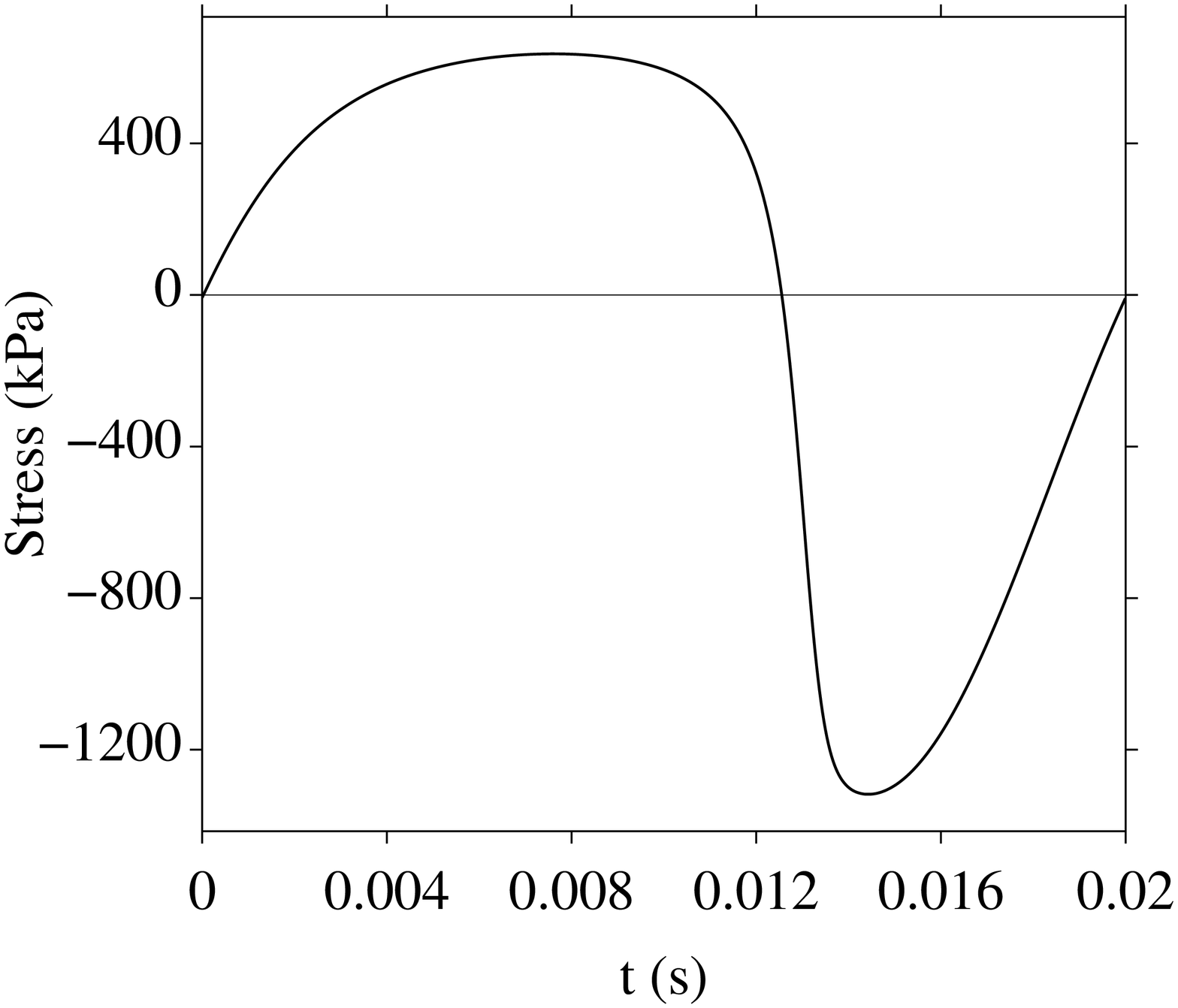}&
\includegraphics[width=6.8cm,height=6cm]{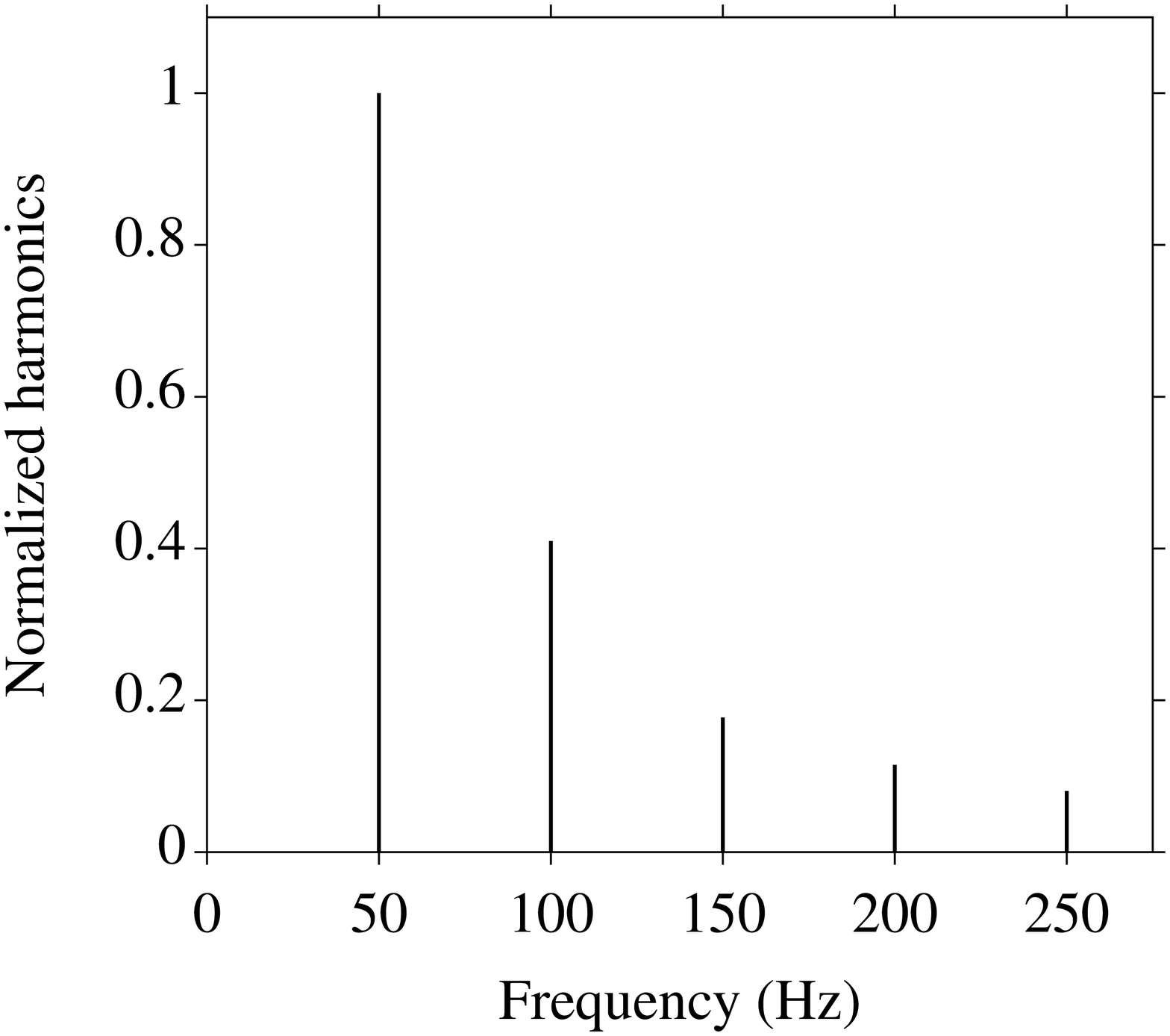}
\end{tabular}
\end{center}
\caption{Test 3: $v_0$ = 0.01 m/s (a), $v_0$ = 0.1 m/s (b), and $v_0$ = 0.4 m/s (c). Left column: numerical values of $\sigma$ during a period; right column: normalized coefficients of the Fourier decomposition.}
\label{Test3}
\end{figure}

As a last test, we simulate a time-harmonic experiment. To do so, we slightly change the system under study: instead of (\ref{IBVP}), we solve
\begin{equation}
\left\{
\begin{array}{l}
\displaystyle
\frac{\textstyle \partial}{\textstyle \partial \,t}\,\boldsymbol{U} + \boldsymbol{A}\, \frac{\textstyle \partial}{\textstyle \partial \,x}\,\boldsymbol{U}=\delta(x-x_s)\,\boldsymbol{S}(t)\quad \mbox{ for } x\in \mathbb{R},\quad x\neq \alpha, \quad t\geq 0,\\
[8pt]
\displaystyle
\boldsymbol{U}(\alpha^+,\,t)=\boldsymbol{D}_0 \left(\boldsymbol{U}(\alpha^-,\,t),\,\frac{\textstyle \partial}{\textstyle \partial\,x}\boldsymbol{U}(\alpha^-,\,t)\right), \\
[8pt]
\displaystyle
\boldsymbol{U}(x,\,0)=\boldsymbol{0}
\quad \mbox{ for } x \in \mathbb{R},
\end{array}
\right.
\label{Source}
\end{equation}
where $\boldsymbol{S}(t)=(0,\,\varepsilon\,\sin(\omega_c\,t))^T$  is a source of mass that acts at $x_s=40$ m (in $\Omega_0$). Except for the source, the parameters are the same than in subsection \ref{SecTest1}. Three values of $\varepsilon$ are considered, yielding the three amplitudes $v_0$: 0.01 m/s (a), 0.1 m/s (b), and 0.4 m/s (c). The periodic values of the transmitted wave are recorded at $x=220$ m; a decomposition into Fourier series is then applied to these values. Figure \ref{Test3} shows the numerical values of $\sigma$ on one period (left column) and the normalized coefficients of the Fourier series (right column). The harmonics generated when $v_0$ increases are clearly seen in this figure.

There are two reasons for displaying these harmonics. First, they help to decide whether it is worthwhile taking the nonlinear effects into account: solving the nonlinear system (\ref{SysNl}) at each time step is more costly than solving a linear system during a pre-processing step, as done with the linear jump conditions (\ref{JClin}) in \cite{ALIMENTAIRE1}. In case (a), the answer is negative; in case (b), it depends on the accuracy required; in case (c), the answer is positive. Secondly, the decrease in the harmonics is linked to the stiffness $K$ and the maximum allowable closure $d$. It may therefore be possible to infer these values by inspecting the harmonics, which is a non-destructive means of evaluating the fracture.

\section{Conclusion}\label{SecConclu}

Here we have studied the propagation of 1-D elastic compressional waves across a contact nonlinearity. The latter feature models a fracture in rocks, but it can also be a useful means of describing other physical situations, such as those encountered in the nondestructive evaluation of material for instance. The present study involves a physical description of the model, a theoretical analysis of its solution, and a numerical time-domain modeling.

The numerical analysis of the interface method needs to be studied further to determine its accuracy and stability (the latter point is still an open question in the linear context). The resolution of the nonlinear system (\ref{SysNl}) needs also to be improved. Lastly, it would be interesting to determine the harmonics generated across the fracture analytically, using the harmonic balance method. 

Many difficulties need to be overcome before tackling realistic 2-D configurations. First, it is required to model the mechanical shear behavior of realistic fractures, along with the compressional behavior \cite{BANDIS83}. These effects can be efficiently modeled in 1-D configurations. Secondly, the 2-D treatment proposed in \cite{ALIMENTAIRE2} to deal with linear contacts is greatly intricated: the linear underdetermined systems of jump conditions are now nonlinear. Thirdly, the computational extra-cost induced by the interface method is substantially higher: the extrapolation matrices used in the interface method need to be computed at each time step and at many grid points along the fracture. Fourthly, if the nonlinear effects are large, a fine grid will be required to converge towards the right solution of the nonlinear systems. Using a local mesh refinement around the fracture might therefore be a useful strategy here.

{\bf Acknowledgments.} We thank Denis Matignon (ENST Paris) and Sergio Bellizzi (LMA Marseille) for their reading of the manuscript and their advice.

\end{document}